\documentclass[10pt]{article}
\usepackage{amsmath}
\usepackage{amsfonts}
\usepackage{amssymb}

\newtheorem{proposition}{Proposition}[section]
\newtheorem{theorem}{Theorem}[section]
\newtheorem{lemma}[theorem]{Lemma}
 
\newtheorem{corollary}[theorem]{Corollary}
\newtheorem{remark}[theorem]{Remark}

\newtheorem{definition}{Definition}

\def\phi{{\varphi}}

\def\B{{\cal B}}
\newcommand{\braket}[2]{\langle #1,#2 \rangle}
\newcommand{\al}{\alpha}

\newcommand{\la}{\lambda}

\newcommand{\ra}{\rightarrow}
\def\phi{{\varphi}}

\DeclareSymbolFont{AMSb}{U}{msb}{m}{n}
\DeclareMathSymbol{\N}{\mathbin}{AMSb}{"4E}
\DeclareMathSymbol{\Z}{\mathbin}{AMSb}{"5A}
\DeclareMathSymbol{\R}{\mathbin}{AMSb}{"52}
\DeclareMathSymbol{\Q}{\mathbin}{AMSb}{"51}
\DeclareMathSymbol{\I}{\mathbin}{AMSb}{"49}
\DeclareMathSymbol{\C}{\mathbin}{AMSb}{"43}

\begin{document}

\title{Schr\"odinger equations and Hamiltonian systems of PDEs with selfdual  boundary conditions}
\author{ Nassif  Ghoussoub\thanks{Partially supported by a grant
from the Natural Sciences and Engineering Research Council of Canada.  } \quad  and \quad Abbas Moameni\thanks{Research supported by a postdoctoral fellowship at the University of British Columbia.}\\
\small Department of Mathematics,
\small University of British Columbia, \\
\small Vancouver BC Canada V6T 1Z2 \\
\small {\tt nassif@math.ubc.ca} \\
\small {\tt moameni@math.ubc.ca}
\\
%\today\\
%\date{January 20, 2005}\\
}
\maketitle

\begin{abstract} Selfdual variational calculus is further refined and used  to address questions of existence of local and global solutions  for various parabolic semi-linear equations, Hamiltonian systems of PDEs, as well as certain nonlinear Schr\"odinger evolutions. This allows for the resolution of such equations  under general  time boundary conditions which include the more traditional ones such as initial value problems, periodic and anti-periodic orbits, but  also yield new ones such as ``periodic orbits up to an isometry" for evolution equations that may not have periodic solutions. In the process, we introduce a method for perturbing selfdual functionals in order to induce coercivity and compactness, while keeping the system selfdual. 

\end{abstract}

\section{Introduction} We develop further the selfdual variational calculus  in order to deal with various parabolic semi-linear equations, Hamiltonian systems of PDEs, as well as certain nonlinear Schr\"odinger evolutions. Our goal  is to solve these equations under general  --sometimes nonlinear-- time boundary conditions which,  besides yielding the more traditional ones such as initial value problems, periodic and anti-periodic orbits, they also yield ``periodic orbits up to an isometry" for certain evolution equations that may not have periodic solutions. We shall use the selfdual variational calculus --developed in \cite{G2, G3, GM3}-- to write these evolution equations  as
\begin{equation}\label{SDE100}
 \left \{ \begin{array}{lcl}
  \hfill  \dot {u}(t)+A u(t)& = & -\bar \partial L (t,u(t)), \\ 
  \hfill  \frac{ u(T) +u(0)}{2}
 & \in & -\bar \partial \ell\big( u(0)-u(T)),
\end{array}\right.
 \end{equation} 
and the Hamiltonian systems as
\begin{eqnarray}\label{SDE200}\left\{ \begin{array}{lcl}
\hfill \dot U(t)+{\cal A}U(t)&= &-J\bar \partial L\big( t,U(t)\big)\\
\hfill \frac {U(T)+U(0)}{2}&\in & -R \bar \partial
\ell \big(U(0)-U(T)\big).
\end{array}\right.
\end{eqnarray}
where  $A$ (resp., ${\cal A}$) is a --non necessarily linear-- operator on a suitable Hilbert space $H$ (resp., $X:=H\times H$),   $J$ is the symplectic operator $J(u,v)=(-v, u)$ and $R$ is the automorphism $R(u,v)=(u, -v)$.  \\ 
The key concept here is the notion of a vector field $\bar \partial L$ that is derived from a convex lower semi-continuous Lagrangian on phase space $L:X\times X^*\to \R\cup\{+\infty\}$ in the following way:  for each $x\in X$, the  --possibly empty-- subset  $\bar \partial L(x)$ of $X^*$ is defined as
 \begin{eqnarray}
\bar \partial L(x): = \{ p \in X^*;  (p, -x)\in \partial L(x,-p) \}.
\end{eqnarray}
 Here $\partial L$ is the subdifferential of the convex function $L$ on $X\times X^*$, which should not be confused with $\bar \partial L$. Of particular interest to us, are those vector fields derived from  {\it anti-selfdual Lagrangians}, i.e., those convex lower semi-continuous Lagrangians $L$ on $X\times X^*$ that satisfy the following duality property:
  \begin{equation}
  L^*(p,x)=L(-x, -p) \quad \hbox{\rm for all $(x,p)\in X\times X^*$},
  \end{equation}
   where here $L^*$ is the Legendre transform in both variables, i.e.,
  \[
  L^*(p,x)= \sup  \{  Re \braket{ y}{\bar p }+  Re  \braket{x}{\bar q }-L(y, q): \,  (y,q)\in X\times X^*\},
  \]
  Such Lagrangians satisfy the following basic property:
 \begin{equation}\label{obs.1}
\hbox{$ L(x, p)+\langle x, p\rangle\geq 0$   for every $(x, p) \in X\times X^{*}$.}
 \end{equation}
  Moreover, 
 \begin{equation}\label{obs.2}
\hbox{  $L(x, p)+\langle x, p\rangle =0$ if and only if $(-p, -x)\in \partial L(x,p)$,}
\end{equation}
which means  that its associated {\it anti-selfdual vector field}  at $x \in X $ is simply
\begin{eqnarray}
\bar \partial L(x):= \{ p \in X^*; L(x,-p)- \langle x,p \rangle=0 \}.
\end{eqnarray}
 Before going further, let us note that  {\it anti-selfdual vector fields} are natural,  but far reaching extensions of subdifferentials of convex lower semi-continuous functions. Indeed, the most
basic anti-selfdual Lagrangians are of the form $L(x,p)= \varphi (x)+\varphi^*(-p)$ where  $\varphi$ is such a function in $X$, and $\varphi^*$
is its Legendre conjugate on $X^*,$ in which case  $\bar \partial  L(x)= \partial \varphi (x).$   More interesting  examples of
anti-selfdual Lagrangians are of the form $L(x,p)= \varphi (x)+\varphi^*(-\Gamma x-p)$  where $\varphi$  is a convex and lower semi-continuous
function on $X,$ and  $\Gamma: X\rightarrow X^*$ is a skew adjoint operator. The corresponding anti-selfdual vector field is then 
$\bar \partial L(x)=\Gamma x+ \partial \varphi (x)$.   Actually, it turned out that every  { \it maximal monotone operator}  (see for example \cite{Br}) is an {\it anti-selfdual vector field} and vice-versa. This fact --proved in \cite{G5}-- means that anti-selfdual Lagrangians can be seen as the {\it potentials} of maximal monotone operators, in the same way as the Dirichlet integral  is the potential of the Laplacian operator (and more generally as any convex lower semi-continuous energy is a potential for its own subdifferential), leading to a variational formulation and resolution of most equations involving maximal monotone operators.   

The main premise of selfdual variational calculus is that many partial differential equations can be formulated as
\begin{equation}
\hbox{$0\in \bar \partial L(x)$ \quad or \quad $-\Lambda x \in \bar \partial L(x)$}
\end{equation}
where $\Lambda:D(\Lambda)\subset X\to X^{*}$ is a linear or non-linear operator, and that solving such an equation amounts to proving that the functional
\begin{equation}
\hbox{$I(x)=L(x,0)$ \quad or \quad $I(x)=L(x,\Lambda x) +\langle x, \Lambda x\rangle$}
\end{equation}
attains its infimum, and --as importantly-- that such an infimum  is equal to  zero. This point of view has been developed in a series of recent papers \cite{G2, G3, GM1, GM2}. However, several new phenomena emerge while dealing with evolutions of the form (\ref{SDE100}) and (\ref{SDE200}), and  many useful new techniques are introduced here  to selfdual variational calculus.  We shall summarize now the main novel ideas, leaving the precise statements and proofs for the following sections. 
  
  \subsection*{(A) The selection of anti-selfdual Lagrangians}
  In applying the general existence results we obtain for equations of the form (\ref{SDE100}) and (\ref{SDE200}), we are often presented with many ways to  associate an anti-selfdual Lagrangian $L$ to the given vector fields. Consider for example, the case of a general  semi-linear evolution equations of the form
\begin{eqnarray}\label{equation.1}
\dot x(t) + A x(t)+w x(t)\in  -\partial \varphi \big( t,x(t)\big)\quad \mbox{ for a.e. }t\in [0,T]
\end{eqnarray}
where 
$w\in \R$, $\varphi (t, \cdot ): H\rightarrow \mathbb{R}\cup \{+\infty\}$ is a proper convex and lower semi-continuous functional on a Hilbert space $H$  and $A: {\rm Dom} (A)\subseteq H \rightarrow H$ is a linear operator. A typical example being  the   complex Ginsburg-Landau equation on $\Omega \subseteq \mathbb{R}^N$, 
\begin{eqnarray} \label{gl.1}
 \begin{array}{lcl}
\frac{\partial  u}{\partial t}-(\kappa +i\alpha)\Delta u + wu= -\partial \phi (t, u(t)) \quad \mbox{ for }t\in (0,T]. 
\end{array}
\end{eqnarray}
We may have several  possible situations: 
\begin{enumerate}
\item {\it The diffusive case}  which corresponds for instance to when $w\geq 0$, $A$ is a positive operator and the --then convex-- function $\Phi (t, x)=\phi (t, x) +\frac{1}{2} \langle Ax, x \rangle +\frac{w}{2}\|x\|_H^2$ is coercive on the right space. In this case, the anti-selfdual Lagrangian is 
\begin{equation}
L(t, x,p)=\Phi (t, x) +\Phi^*(t, -A^ax-p)
\end{equation}
where $A^a$ is the anti-symmetric part of the operator $A$. 
\item {\it The non-diffusive case} which essentially means that one of the above requirements is not satisfied, e.g., $w<0$ or if $A$ is unbounded and purely skew adjoint ($\kappa=0$).  The anti-selfdual Lagrangian is then 
\begin{equation}
L(t, x,p)=e^{-2\omega t}\left\{\phi(t, e^{\omega t}S_tx)+ \phi^*(t, -e^{\omega t}S_tp)\right\}
\end{equation}
where $S_t$ is the $C_0$-unitary group associated to the skew-adjoint operator $A$. 
This non-diffusive case cannot be formulated on ``energy spaces" and therefore requires less stringent  coercivity conditions. However, the equation may not  in this case have solutions satisfying the standard boundary conditions. Instead,  and as we shall see below, one has to settle for solutions that are periodic but only up to the isometry $e^{-TA}$.
\item {\it The mixed case} which deals with 
\begin{eqnarray}\label{equation.1}
\dot x(t) + A_1 x(t)+A_2 x(t) +w x(t)\in  -\partial \varphi \big( t,x(t)\big)\quad \mbox{ for a.e. }t\in [0,T]
\end{eqnarray}
where $A_1$ is a bounded positive operator  and $A_2$ is an unbounded and purely skew adjoint operator.  One example we consider,  is  the  following evolution equation with an advection term.
\begin{eqnarray}
\dot u (t)+a.\nabla u(t)-i\triangle u+ w u(t)=-\partial \varphi (t, u(t)) \quad \mbox{ for } t\in [0,T].
\end{eqnarray}
The anti-selfdual Lagrangian is then 
\begin{equation}
L(t, x,p)=e^{-2\omega t}\left\{\phi(t, e^{\omega t}S_tx)+ \phi^*(t, -e^{\omega t}A_1^aS_tx-e^{\omega t}S_tp)\right\}
\end{equation}
where $S_t$ is the $C_0$-unitary group associated to the skew-adjoint operator $A_2$. Again, one then gets the required boundary condition  up to the isometry $e^{-TA_2}$.

 \end{enumerate} 
 \subsection*{(B) The selection of boundary Lagrangians}   The interior Lagrangians $L$  above --and throughout this paper--  are expected in the applications to be smooth  and hence their subdifferentials will coincide with their  differentials, and the corresponding inclusions will  often be equations.  It is however crucial here 
that
the boundary Lagrangians $\ell$ be  allowed to be degenerate so that they can cover the boundary conditions that we now discuss. Indeed, the selfdual boundary conditions in  (\ref{SDE100}) 
often translates into
\begin{equation} \label{SDB}
\hbox{$\frac{ v(0)+e^{-wT} S_{-T} v(T)}{2} \in -\partial \psi \left(v(0)- e^{-wT} S_{-T} v(T) \right)$} 
\end{equation} 
where $\psi$
is a convex function on $H$,  and $(S_t)_t$ is the $C_0$-unitary group associated to the skew-adjoint part of the operator. Here is a sample of the  various boundary conditions that one can obtain by choosing $\psi$ accordingly in (\ref{SDB}). 
 \begin{enumerate}
\item  {\it Initial boundary condition,} say $v(0)=v_0$ for a given $v_0\in
H$, then it suffices to choose $\psi (u)= \frac{1}{4}\|u\|_H^2-\braket{u}{v_0}$.
 
\item {\it Periodic type solutions} of the form $ v(0)= S_{-T} e^{-wT} v(T)$, then $\psi$ is chosen as:
\begin{eqnarray*}\psi(u)=\left\{\begin{array}{ll}
0 \quad &u=0\\
+\infty &\mbox{elsewhere}.\end{array}\right.
\end{eqnarray*}

\item {\it Anti-periodic type solutions} $ v(0)= -S_{-T} e^{-wT} v(T)$,  then $\psi (u)=0$ for each $u \in H.$

\end{enumerate}
 In the latter cases, we shall say that the solutions are periodic and anti-periodic orbits up to an isometry.  

 \subsection*{(C) The use of selfduality to induce coercivity and compactness}
 Typical Hamiltonian systems of PDEs we are aiming to solve via a selfdual variational approach are:
\begin{equation}\label{ham-ex11}
 \left\{ \begin{array}{lcl}
-\dot v(t)-\Delta (v+u)+b.\nabla v&= & \partial \phi_1 (t,u)   \\
\hfill \dot u(t)-\Delta (u+v)+ a.\nabla u&= &\partial \phi_2(t,v) 
\end{array}\right. 
\end{equation}
as well as 
\begin{eqnarray}\label{ham-ex21}\left\{ \begin{array}{lcl}
-\dot v(t)+\Delta^2 v - \Delta v&= & \partial \phi_1 (t,u)  \\
\hfill  \dot u(t)+\Delta^2 u+ \Delta u&= &\partial \phi_2(t,v)  
\end{array}\right.
\end{eqnarray}
with Navier-type state boundary conditions, and where $\phi_i, i=1,2$ are convex functions on some $L^p$-space. Now, in order to deal with such systems, one needs to overcome the fact that the cross product $u\to \int_0^T\langle u(t), \dot u (t)\rangle \, dt$ is not necessarily weakly continuous as in the case of finite dimensional Hamiltonian systems. One important novelty in this paper, is the introduction of a way  to perturb a selfdual functional so as to make it coercive in an appropriate space without destroying selfduality.  We shall now illustrate the main ideas on the following simplified example:
\begin{equation}\label{sample}
 \Gamma x +Ax = -\partial \phi (x) 
\end{equation}
where $\phi$ is a convex lower semi-continuous function on a Hilbert space $H$, and where $A:D(A)\subset H \to H$ and $\Gamma:D(\Gamma)\subset H \to H$ are linear operators. The most basic selfdual functional associated to (\ref{sample}) is
\begin{equation*}
I(x)=\phi (x) +\phi^*(-Ax-\Gamma x)+\langle x, Ax +\Gamma x\rangle. 
\end{equation*}
The main ingredients that allow to show that the infimum is zero and that it is attained,  are:
\begin{enumerate}
\item The weak lower semi-continuity of the function $x\to \langle x, Ax +\Gamma x\rangle$ on $D(A)\cap D(\Gamma)$, and
\item A coercivity condition which implies for example that  $\lim_{\|x\|\to +\infty} I(x) =+\infty$. 
%$\lim_{\|x\|\to +\infty}\phi^*(-Ax-\Gamma x)+\langle x, Ax +\Gamma x\rangle =+\infty$. 
\end{enumerate}
Now suppose that $A$ satisfies
%\begin{equation*}\hbox{
$\langle Ax,x\rangle \geq c_0 \|x\|^2$ for all $x\in D(A)$, 
%}\end{equation*}
and that $A^{-1}$ is a compact operator,  then one can strengthen the topology on the domain of the functional $I$ by considering the Hilbert space $Y$ that is the completion of $D(A)$ for the norm $\|u\|_A=\langle Au, u \rangle$ induced by the scalar product  $\langle u, v \rangle_Y=\langle u, Av \rangle_H$.
  Note since the injection of $Y$ into $H$ is  compact, the map $x\to \langle  x, Ax\rangle$ is readily weakly  continuous on $Y$,  and  the function $x\to \langle x,  \Gamma  x \rangle$ has a better chance to be lower semi-continuous for the weak topology of $Y$. 
 On the other hand, by considering $I$ on the space $Y$, we often lose coercivity for the new norm, which is not guaranteed by the following sub-quadratic growth that we assume on $\phi$.    
\begin{equation}
\hbox{$-C\leq \phi (u) \leq \frac{\beta }{2}(\| u \|^2 +1)$ for   $u\in H$,}
\end{equation} 
 for some $\beta>0$ and $C\in \R$. Indeed, such a condition yields 
  \[
\phi^*(-Ax  - \Gamma x) \geq \frac{1}{2\beta}(\| Ax +\Gamma x\|^2 -1) \geq  \frac{1}{2\beta}\| Ax\|^2 + \frac{1}{\beta}\langle Ax, \Gamma x \rangle -\alpha, 
\]
in such a way that the functional $I(x)- \frac{1}{\beta} \langle Ax, \Gamma x \rangle$ is coercive for the norm of $Y$. But this new functional is however not selfdual, and so to remedy this, we use the fact that often, the cross product $\langle Ax, \Gamma x\rangle$ can be resolved via a Green-Stokes type formula of the form:
\begin{equation}
\hbox{$ \langle Ax, \Gamma x \rangle +\langle T\B x, R \B x \rangle =0$ for all $x\in Y$},
\end{equation} 
where $\B: D(\B)\subset H\to H_0$ is an  operator into a boundary Banach space $H_0$, $T$ is an operator on $H_0$ and  $R:H_0\to H_0^*$ is such that for some $c>0$, 
\begin{equation}
\hbox{$|\langle \B x, R \B x \rangle| \leq c\|x\|_Y^2$ for all $y\in Y$.}
\end{equation}
We then consider any convex lower semi-continuous function $\psi$ on $H_0$, and let $\ell (a,b)=\psi (a) +\psi^*(-Tb)$ for $(a, b) \in H_0\times H_0^*$ in such a way that 
\begin{equation}
\hbox{$\ell (a,b)\geq -\langle Ta, b\rangle$ for all $(a,b)\in H_0\times H_0^*$.}
\end{equation} 
The following   functional
\[
J(x)=I(x) + \frac{1}{\beta}\ell (\B x, R\B x)+(\frac{1}{\beta}\ell)^*(-R\B x, -\B x)+2\langle \B x, R \B x \rangle
\]
is then non-negative, selfdual,  but also coercive on $Y$ as soon as $\beta <\frac{c_0}{2c}$ since 
\[
J(x)\geq I(x) +\frac{c_0}{\beta}\|x\|_Y^2  -2 c\|x\|_Y^2- C.
\]
The infimum of  $J$ on $Y$ is then equal to zero and is attained  at a point $u\in Y$ satisfying
 \begin{eqnarray}
\label{cute}
 \left\{ \begin{array}{lcl}
Au +\Gamma u &=& -\partial \phi (u)\\
\hfill R\B u&\in& \frac{-1}{\beta} \partial \psi (\B u). 
\end{array}\right.
\end{eqnarray}
It is worth noting that the required bound $\beta <\frac{c_0}{2c}$ normally leads to a time restriction in evolution equations and often translates into local existence results as opposed to the global ones in the case of $(\ref{SDE100})$.
 The relevance of this approach will be illustrated in the section on Hamiltonian systems of PDEs. \\
 \subsection*{(D) Schr\"odinger evolutions and nonlinear selfdual principles:}
  In the case of a  Schr\"odinger equation of the form 
   \begin{eqnarray} \label{sch.1}
  \begin{array}{lcl}
i\frac{\partial  u}{\partial t}-(\kappa +i\alpha)\Delta u  + wu= \partial \phi (t, u(t)) \quad \mbox{ for  }t\in (0,T] \\
 \end{array}
\end{eqnarray}
 where  $w\in \R$, $\kappa$ and  $\alpha \geq 0$, and $\varphi (t, \cdot)$
is a proper convex and lower semi-continuous functional on  $L^2(\Omega)$ or $H_0^1(\Omega)$,  we shall rewrite it  in the form 
 \begin{eqnarray} \label{sch.2}
  \begin{array}{lcl}
-\frac{\partial  u}{\partial t} - A u(t) -\Lambda u (t)=  \partial \Psi (t, u(t)) \quad \mbox{ for  }t\in [0,T] 
\end{array}
\end{eqnarray}
where $A:=-i\kappa \Delta$ is a skew adjoint operator, $\Lambda:=iwu-i\partial \phi (t, u(t)) $ is a nonlinear operator, while $\Psi (u)=\frac{\alpha}{2} \int_\Omega |\nabla u|^2 dx$ . Here again, there are two ways for ``embedding"  the skew-adjoint operator $A$ into an anti-selfdual Lagrangian, so as to reduce it to a nonlinear evolution of the form 
\begin{equation}\label{SDE1}
 \begin{array}{lcl}
  \hfill  \dot {u}(t)+\Lambda u(t)& = & -\bar \partial L (t,u(t)) 
\end{array}
 \end{equation} 
where $\Lambda$ is a nonlinear operator. This latter equation was dealt with in \cite{GM3} in the context of the Navier-Stokes evolutions, but we show here how it can be combined with semi-group theory in order to handle nonlinear evolutions with an additional  skew-adjoint term.

The paper is organized as follows. We start by reviewing in section 2, some basic properties of seldual Lagrangians and functionals.  In section 3,  we establish a selfdual variational principle for semi-linear parabolic equations with general boundary conditions. Applications to complex Ginsburg-Landau evolutions,   coupled flows and other wave-type equations are given. Section 4 is concerned with  Hamiltonian systems of PDEs, where additional selfdual terms are used to induce coercivity and compactness,  while section 5 deals with nonlinear evolutions and in particular Schr\"odinger equations. Most of this paper is self-contained, though it is preferable to read it in conjunction with \cite{G2}, \cite{G3} which introduce the basics about selfduality and its immediate applications. Section 5 is however heavily dependent on \cite{GM3}.
 
\section{Basic properties of selfdual functionals} 

We start by recalling the concept of an anti-selfdual Lagrangian and its main properties.   Let $X$ be a (real or complex) reflexive Banach space and let $X^*$ be its dual. Hence forth, we shall simply denote the real scalar product  $Re\braket{\ }{\ }$    by $\braket{\ }{\ }.$

Given a function on phase space $L:X\times X^*\to \R\cup\{+\infty\}$, we define the  {\it derived vector field of } $L$ at $x \in X $ to be the -possibly empty- subset  of $X^*$ given by:
\[
\bar \partial L(x) =\{p\in X^*; L(x,-p)+L^*(p, -x)=2\langle x, p\rangle\}.
\]
If $L$ is convex and lower semi-continuous on $X\times X^*$, then 
\[
\bar \partial L(x) =\{p\in X^*;  (p, -x)\in \partial L(x,-p)\}.
\]
 If now $L$ is an anti-selfdual Lagrangian, then 
 \[
\bar \partial L(x) =\{p\in X^*;  L(x,-p)-\langle x, p\rangle =0\}.
\]
The  Hamiltonian (resp. co-Hamiltonian) $H=H_L$ on $X\times X$ (resp. $\tilde H=\tilde H_L$ on $X^*\times X^*$) corresponding to $L$  are given by: 
\begin{eqnarray*}
\hbox{$H_L (x,y)=\sup \{ \langle y,p \rangle-L(x, p); p\in X^*\}$ \quad (resp., \quad 
 $ \tilde H_L (p,q)=\sup \{ \langle x,q \rangle-L(x, p); x\in X\})$}
\end{eqnarray*}
 
\subsection*{Basic variational principles for selfdual functionals}

Our main premise is that many partial differential equations can be formulated as
\begin{equation}
\hbox{$0\in \bar \partial L(x)$ \quad or \quad $\Lambda x \in -\bar \partial L(x)$}
\end{equation}
where $\Lambda:D(\Lambda)\subset X\to X^{*}$ is a linear or non-linear operator, and that solving such an equation amounts to proving that the functional
\begin{equation}
\hbox{$I(x)=L(x,0)$ \quad or \quad $I(x)=L(x,\Lambda x) +\langle x, \Lambda x\rangle$}
\end{equation}
attains its infimum, and --as importantly-- that such an infimum  is equal to  zero. 

\begin{definition}\rm A functional $I: X\to \R\cup \{+\infty\}$ is said to be {\it completely selfdual on $X$},  if there exists an anti-selfdual Lagrangian $L$ on $X\times X^*$ such that $I(x)=L(x,0)$ for every $x\in X$. 
\end{definition}
Note that  completely selfdual functionals can also be written as 
  \begin{equation} \label{sd}
 \hbox{$I(x)=\sup\limits_{y\in X}H_L(y, -x)$ \quad for all $x\in X$,} 
 \end{equation}
where $H_L$ is the Hamiltonian associated of $L$. 
  The function $M(x,y)=H_L(y, -x)$  has some remarkable properties. In particular, it satisfies:
  \begin{enumerate}  
  \item For each $y\in X$, the function $x\to M(x,y)$  is weakly lower semi-continuous;   
\item  For each $x\in X$, the function $y\to M(x,y)$ is concave;
\item For each $x \in X$, we have $M(x,x)\leq 0$. 
\end{enumerate}
Such an $M$ will be called an {\it anti-symmetric Hamiltonian} on $X\times X$.
\begin{definition} \rm We say that a functional  $I: X\to \R^+ \cup \{+\infty\}$  on a  Banach space $X$ 
  is {\it selfdual on a convex set $D\subset X$},   if there exists an anti-symmetric  Hamiltonian $M: D\times D \to \R$ such that 
\begin{equation}
\hbox{$I(x)=\sup\limits_{y\in D}M(x,y)$ for every $x\in D$.}
\end{equation}
\end{definition}
The following two existence results will be frequently used in the sequel. They give sufficient conditions for the infimum of selfdual functionals  to be attained, and --as importantly-- to be zero.  
 
\begin{theorem} {\rm \cite{G2}} \label{one} Let $I$ be a completely selfdual functional on a reflexive Banach space $X$, such that its associated anti-selfdual Lagrangian $L$ on $X\times X^{*}$ satisfies for some $x_{0}\in X$, that   $p\to L(x_{0},p)$ is bounded above on  a neighborhood of the origin in $X^{*}$. 
  Then  there exists $\bar x\in X$ such that   $
I(\bar x)=\inf_{x\in X}I(x)=0.$ 
  \end{theorem}

\begin{theorem} \label{two} {\rm \cite{G3} } Let  $I$ be a  selfdual functional  on a convex  closed subset $D$ of a reflexive Banach space $X$, such that its associated anti-symmetric Hamiltonian $M$ on $D\times D$  satisfies $\lim\limits_{\|x\|\to \infty} M(x, x_0)=+\infty$ for some $x_0\in D$.   Then there exists $\bar x\in D$ such that $
I(\bar x)=\inf_{x\in D}I(x)=0.$
\end{theorem} 

\subsection*{Operations on selfdual Lagrangians} 

 We now summarize various permanence properties enjoyed by the class of anti-selfdual Lagrangians.   For the proofs, we  refer to \cite{G2}.
\begin{proposition} \label{permanence}
Suppose  $L$ is an {\it anti-selfdual Lagrangian} on
$ X\times X^*$, where $X$ is a reflexive Banach space, then
\begin{enumerate}
 \item
  For every $\mu >0$, the Lagrangian ${L_\mu}(u,p):=\mu^{-2}
L(\mu u, \mu p)$  is also  anti-selfdual. 
 \item If $A:X\to X^*$ is a bounded skew adjoint operator, then the Lagrangian
 $M(u, p)=L(u, Au+p)$ is again {\it anti-selfdual
 Lagrangian}.
  \item If  $X=H=X^*$ and $S$ is a  unitary operator on  a Hilbert space $H$ (i.e. $SS^*=S^*S=I$), then the Lagrangian $
L_S(u,p):= L(Su,Sp)$  is also an {\it anti-selfdual
 Lagrangian}.

 \end{enumerate}
\end{proposition}
 Suppose now that we have an evolution triple $X\subset H\subset X^*$, where $X$ is reflexive, $H$ is a Hilbert space and where each space is dense in  the following one.  Also assume that there exists  a  linear and symmetric duality map $D$ between $X$ and $X^*$, in such a way that $\|x\|^2=\langle x,  Dx \rangle$.  We can then consider $X$ and $X^*$ as Hilbert spaces with the following inner products,
\begin{eqnarray}
\langle u,v \rangle_{X \times X}:=\langle D u,v \rangle \quad \text{ and } \quad 
\langle u,v \rangle_{X^* \times X^*}:=\langle D^{-1}u,v \rangle
\end{eqnarray}
A typical example is the evolution triple  $X=H_0^1(\Omega) \subset H:= L^2(\Omega)\subset X^*=H^{-1}(\Omega)$ where the duality map is given by $D= - \triangle.$  \\
If now $\bar S$ is an isometry on $X^*$, then $S=D^{-1}\bar S D$ is also an isometry on $X$, in such a way that 
\begin{equation}
\hbox{$\langle u,p \rangle=\langle S_t u , \bar S_t p\rangle$ for all $u \in X$ and $p \in X^*$.}
\end{equation}
Indeed, we have
\begin{eqnarray*}
\langle S u , \bar S p\rangle=\langle DS u , \bar S p\rangle_{X^* \times X^*}=\langle  \bar S D u , \bar S p\rangle_{X^* \times X^*}=\langle  D u ,  p\rangle_{X^* \times X^*}=\langle u,p \rangle.
\end{eqnarray*}
from which we can deduce that 
\begin{eqnarray*}
\| S u\|_X^2= \langle S u, S u \rangle_{X \times X}= \langle S u, DS u \rangle= \langle S  u, \bar S Du \rangle= \langle u, Du \rangle=\|u\|_X^2.
\end{eqnarray*}
Moreover, if $L$ is an anti-selfdual Lagrangian on $X\times X^*$, then  $L_{S}:=L(S u, \bar S p)$ is also an anti-selfdual Lagrangian on $X \times X^*$, since
\begin{eqnarray*}
L^*_{S}(p,u)&=&\sup\{\langle v,p \rangle+ \langle u, q \rangle -L_{S}(v,q); (v,q) \in X \times X^*\}\\
&=&\sup\{\langle S v, \bar S p \rangle+ \langle S u, \bar S q \rangle -L(S v,\bar S q); (v,q) \in X \times X^*\}\\
&=&L^*(\bar S p,S u)=L(-Su,-\bar S p)=L_{S}(-u,-p).
\end{eqnarray*}
 We shall  also make repeated use of the following  lemma which describes three ways of regularizing an anti-selfdual Lagrangian by inf-convolution. It is an immediate consequence of the calculus of anti-selfdual Lagrangians developed in \cite{G2} to which we refer the reader. 
 \begin{lemma} \label{2-reg} For a  Lagrangian  $L: X \times X^*\rightarrow \mathbb{R}\cup \{+\infty\}$, define for every $(x, r)\in X\times X^*$
\begin{eqnarray*}
L^1_{\lambda}(x,r) =\inf \{ L(y,r)+ \frac {\|x-y\|^2}{2\lambda }+ \frac {\lambda \|r\|_*^2}{ 2}; y \in X \}
\end{eqnarray*}
and
\begin{eqnarray*}
L^2_{\lambda}(x,r) =\inf \{ L(x,s)+ \frac {\|r-s\|_*^2}{2\lambda }+ \frac {\lambda \|x\|^2}{2}; s \in X^* \}
\end{eqnarray*}
and 
\begin{eqnarray*}
L^{1,2}_{\lambda}(x,r)=\inf \big \{   L(y,s)+ \frac{1}{2 \la}\|x-y\|^2+  \frac{\la}{2 }\|r\|_*^2+ \frac{1}{2 \la}\|s-r\|_*^2+  \frac{\la}{2 }\|y\|^2; \, y \in X, s \in X^*\big \}
\end{eqnarray*}

If $L$ is anti-selfdual then the following hold:
\begin{enumerate}
\item   $L^1_{\lambda}$, $ L^2_{\lambda}$ and $L^{1,2}_{\lambda}$ are  also anti-selfdual Lagrangians on $X \times X^*$. 
\item  $L^1_{\lambda}$ (resp., $L^2_{\lambda}$) (resp., $L^{1,2}_{\lambda}$) is continuous in the first variable (resp., in the second variable) (resp., in both variables).
\item  $H_{L^1_{\lambda}}$ and $\tilde H_{L^2_{\lambda}}$ are continuous in both variables. 
\item Suppose $L$ is bounded from below.  If $x_{\la}\rightharpoonup x$ and $p_{\la}\rightharpoonup p$ weakly in $X$ and $X^*$ respectively as $\lambda \to 0$,  and if  $L^{1,2}_{\lambda}(x_{\la},p_{\la})$ ( resp.,  $L^{1}_{\lambda}(x_{\la},p_{\la})$)  ( resp., $L^{2}_{\lambda}(x_{\la},p_{\la})$ is bounded from above,  then
\begin{eqnarray*}
L(x,p) & \leq &  \liminf_{\la\to 0}L^{1,2}_{\lambda}(x_{\la},p_{\la})\\  resp.,  \quad L(x,p) & \leq &\liminf_{\la\to 0}L^{1}_{\lambda}(x_{\la},p_{\la})\\   resp.,  \quad L(x,p) & \leq & \liminf_{\la\to 0}L^{2}_{\lambda}(x_{\la},p_{\la}).
\end{eqnarray*}

\end{enumerate}
\end{lemma}
We shall  make frequent use of  the following  lemma. 

\begin{lemma} \label{mild.estimate} Let $X \subseteq H \subseteq X^*$ be an evolution triple and let $L$ be  an anti-selfdual Lagrangian on $X \times X^*.$
\begin{enumerate}
\item  Assume that for $C>0$  and $r>1 $, we have
$-C\leq L (x,0)\leq C (1+ \|x\|_X^r)$ for all $x\in X$,
 then there exist $C_1>0$ and $C_2>0$ such that
$L(x,q )\geq C_1\|q\|_{X^*}^s- C_2$  
for every $(x,q ) \in X \times X^*$, where  $\frac
{1}{r}+\frac {1}{s}=1$.
\item  Assume that for $C_1, C_2>0$  and $r_1 \geq r_2>1 $ we have
$C_1 ( \|x\|_X^{r_2}-1)\leq L (x,0)\leq C_2 (1+ \|x\|_X^{r_1})$ for all $x\in X$,
then $L$ is continuous in both variables and
the following Lagrangian
\begin{eqnarray*}
M(u,p):=
\left\{ \begin{array}{lcl}
L(u,p), \quad \quad \quad u\in X,\\
+\infty \quad \quad \quad u\in H \setminus X,
\end{array}\right.
\end{eqnarray*}
is anti-selfdual on $H \times H.$
\end{enumerate}
\end{lemma}

\paragraph{Proof:}  (1) For $(x,q) \in X \times X^*$ we have,
\begin{eqnarray*}
L(x,q)&=&\sup_{(y,p) \in X \times X^*} \left\{ \braket {x}{p}+\braket {y}{q}-L^*(p,y)\right\}\\
         &=&\sup_{(y,p) \in X \times X^*} \left\{ \braket {x}{p}+\braket {y}{q}-L(-y,-p)\right\}\\
 &\geq & \sup_{y \in X} \left\{\braket {y}{q}-L(-y,0)\right\}\\
 &\geq & \sup_{y \in X} \left\{\braket {y}{q}-C (1+ \|y\|_X^r)\right\}\\
&=&C_1\|q\|_{X^*}^s- C_2
\end{eqnarray*}
for some positive constants $C_1$ and   $C_2.$

To prove part (2),  note  first that the given coercivity and
bounded assumptions on $L(x,0)$ ensures the boundedness of $L(.,.)$
in $X \times X^*$ and therefore the continuity. Indeed, for some
$C_1, C_2>0$ we have
\begin{eqnarray*}
C_1 ( \|p\|_X^{s_1}+\|x\|_X^{r_2}-1)\leq L (x,p)\leq C_2 (1+ \|x\|_X^{r_1}+ \|p\|_X^{s_2}) \quad \quad \quad (\frac
{1}{r_i}+\frac {1}{s_i}=1, \quad i=1,2).
\end{eqnarray*}
 Now we prove that $M$ is an  anti-self dual Lagrangian on $H \times H.$ Indeed,  fix $(x,q) \in H \times H.$
If $x \in X$, then
\begin{eqnarray*}
M^*(-q,-x)&=&\sup_{(y,p) \in H \times H} \left\{ \braket {-x}{p}+\braket {y}{-q}-M(y,p)\right\}\\
&=&\sup_{(y,p) \in X \times H} \left\{ \braket {-x}{p}+\braket {y}{-q}-L(y,p)\right\}.
\end{eqnarray*}
Since $L(\cdot, \cdot)$ is continuous and $H$
 is dense in $X^*$, we have
\begin{eqnarray*}
  M^*(-q,-x)  &=&\sup_{(y,p) \in X \times X^*} \left\{ \braket {-x}{p}+\braket {y}{-q}-L(y,p)\right\}\\
&=& L^*(-q,-x)=L(x,q)=M(x,q).
\end{eqnarray*}
Now if $x \not \in X$, then
\begin{eqnarray*}
M^*(-q,-x)&=&\sup_{(y,p) \in H \times H} \left\{ \braket {-x}{p}+\braket {y}{-q}-M(y,p)\right\}\\
&=&\sup_{(y,p) \in X \times H} \left\{ \braket {-x}{p}+\braket {y}{-q}-L(y,p)\right\}\\
    &\geq & \sup_{p \in H} \left\{\braket {-x}{p}-L(0,p)\right\}\\
 &\geq & \sup_{p \in H} \left\{\braket {-x}{p}-C (1+ \|p\|_{X^*}^{s_2})\right\}\\
&=&+\infty =M(x,q).
\end{eqnarray*}

\subsection*{Time-dependent selfdual Lagrangians}

\begin{definition}  A  time dependent Lagrangian on $[0,T]\times X\times X^*$ is any
   function $L: [0,T]\times X\times X^*\to \R \cup \{+\infty\}$
   that is  measurable with respect to   the  $\sigma$-field  generated by the products of Lebesgue sets in  $[0,T]$ and Borel sets in $X\times X^*$.
 We shall say that  such a Lagrangian $L$ is  {\it anti-selfdual}  on $[0,T]\times X\times X^*$  if  for any $t\in [0,T]$, the map $L_t:(x,p)\to L(t, x,p)$ is an anti-selfdual Lagrangian on $X\times X^*$.
\end{definition} 
Let $H$ be a  Hilbert space with $\braket{\ }{\ }$ as scalar product over a real or a complex field.
Let $[0,T]$ be a fixed real interval and consider the space $L^{2}_{H}$
of  integrable functions from $[0,T]$ into $H$ with norm
    $ \|u\|^2_{L^2_H} =(\int_0^T \|u(t)\|_{H}^{2}
    dt)^{\frac{1}{2}}.$
 Consider the Hilbert path space
$A_H^2=\left\{ u:[0,T]\ra H; \dot u\in L_H^2\right\}$
consisting of all absolutely continuous arcs $u:[0,T]\ra H$,
equipped with the norm
\[
\| u\|_{A_H^2}=\left( {\left\|\frac{u(0)+u(T)}{2}\right\|}_H^2
+\int_0^T \|\dot u\|^2\, dt\right)^{\frac{1}{2}}.
\]
We shall  identify $A_H^2$ with the product space $H\times L_H^2$, in such a way that its dual $(A_H^2)^*$ can also be identified
 with $H\times L_H^2$ via the formula
\begin{eqnarray*}
{\braket{u}{(p_1,p_0)}}_{A_H^2,H\times L_H^2}
  =Re  \braket{\frac{u(0)+u(T)}{2}}{p_1} +\int_0^T
  Re  \braket{\dot u(t)}{p_0(t)}\, dt
\end{eqnarray*}
where $u\in A_H^2$ and $(p_1,p_0)\in H\times L_H^2$. The following was proved in \cite{GM1}.

\begin{proposition}\label{lift}
Suppose $L$ is an anti-selfdual Lagrangian on $[0,T]\times H\times H$ and that $\ell$ is an anti-selfdual
Lagrangian on $H\times H$, then the Lagrangian defined on $A_H^2\times {(A_H^2)}^*=A_H^2\times (H\times L_H^2)$ by
\begin{eqnarray*}
{\cal M}(u,p)=\int_0^T L\big( t,u(t)+p_0(t),\dot u(t)\big)\, dt
  +\ell \left( u(0)-u(T)+p_1,\frac{u(0)+u(T)}{2}\right)
\end{eqnarray*}
is anti-selfdual Lagrangian on $A_H^2\times (L_H^2\times H)$.
\end{proposition}
We shall need the following facts about semi-groups of operators. 
  \begin{definition}
A $C_0$-group on $H$ is a family of bounded operators
$S=\{S_t\}_{t\in \mathbb{R}}$ satisfying
\begin{description}
\item[$(i)$]  $S_tS_s=S_{t+s}$ for each $t,s \in \mathbb{R},$
\item[$(ii)$]  $S(0)=I,$
\item[$(iii)$]  The function $t\to S_tu \in C(\mathbb{R}, H)$ for each $u \in H.$
\end{description}
  \end{definition}
 We recall a celebrated result of Stone.  
\begin{proposition} An operator $A: D(A)\subset H\to H$ on a Hilbert space $H$ is skew-adjoint if and only if it is the infinitesimal generator of a $C_0$-group of unitary operators $(S_t)_{t\in \R}$ on $H$. In other words, we have $A x=\lim\limits_{t\downarrow 0} \frac{S_tx-x}{t}$ for every $x\in D(A)$. 
\end{proposition}
We shall sometimes denote the  group $S_t$ by $e^{tA}$. It follows from the above that if  $(S_t)_t$ is such a group and if $L$ is a time dependent anti-selfdual Lagrangian on $[0,T]\times H\times H$, then so is the Lagrangian $L_S(t, u, p):= L(t, S_tu, S_tp)$.
 
The same holds if $X\subset H\subset X^*$ is an evolution triple with a linear and symmetric duality map $D$. Indeed, let $(\bar S_t)_{t \in \R}$ be a $C_0-$unitary group of operators associated to a skew-adjoint operator $A$ on the dual space $X^*$ viewed as a Hilbert space (with scalar product $\langle D^{-1}p, q\rangle$).  By defining the  maps $(S_t)_{t \in \R}$ on $X$ via the formual $S_t=D^{-1}\bar S_tD$, we deduce from the above that if $L$ is a time dependent anti-selfdual Lagrangian on $[0,T]\times X\times X^*$, then so is the Lagrangian $L_S(t, u, p):= L(t, S_tu, \bar S_tp)$.
 
\section{Selfdual variational principles for parabolic equations}
This section is concerned with existence results for evolutions of the form
\begin{eqnarray}
 \left\{ \begin{array}{lcl} \hfill \dot{u}(t) &=&- \bar \partial
   L\big(t,u(t))\quad\forall t\in [0,T]\label{E200}\\
 \frac{u(0)+u(T)}{2}  & \in& -\bar \partial \ell (u(0)-u(T)).\label{E300}
 \end{array}\right.
\end{eqnarray}
where $L$ and $\ell$ are anti-selfdual Lagrangians. We then apply it to equations of the form 
\begin{eqnarray}
 \label{BC200}
\dot{ u}(t) + A u(t) +\omega u(t) &=&-\partial\phi( t, u(t))\quad \mbox{ for a.e. }t\in [0,T] \\
\frac{ u(0)+e^{TA}e^{-wT} u(T)}{2}
   &\in& - \partial \psi \left( u(0)- e^{TA}e^{-wT} u(T)\right),   
\end{eqnarray}
where $\phi$ and $\psi$ are convex functions, $A$ is a skew-adjoint operator and $w\in \R$. 
Such principles were developed  in \cite{G2} and \cite{GT1} for initial-value problems associated to 
(\ref{BC200}), while more general boundary conditions
were dealt with in \cite{GM2} but only in the case of a gradient flow (i.e., when $A=0$).

We start with the following proposition. 
\begin{proposition} \label{prop.1}Suppose  $L$ is a time-dependent anti-selfdual Lagrangian  on $[0,T]\times H\times
H$ and  let $\ell $ be an   anti-selfdual Lagrangian on  $H\times H.$
Assume the following conditions:
 \begin{description}
\item[($A'_1$):] For some $n>1$ and $C>0$ we have
 $-C<\int_0^T L(t,x(t),0)\, dt\leq C\big(\| x\|_{L_H^2}^n+1\big)$ for all $x\in L_H^2$.
\item[$(A'2)$:]  $\int_0^T L\big( t, x(t),p(t)\big) \,
dt\rightarrow\infty$ as $\|x\|_{L^2_H} \rightarrow\infty$ for every
$p \in L^2_H.$
\item[$(A'3)$:] $\ell$  is bounded from below and $0 \in {\rm Dom}(\ell).$
\end{description}
Then the functional $ I(x)=\int_0^T L\big( t, x(t),\dot{x}(t)\big)\,
dt
  +\ell \big(x(0)-x(T),
  \frac{x(0)+ x(T)}{2}\big)
$
 attains its minimum at a path $u\in A_H^2$ satisfying
\begin{eqnarray}
I(u)&=&\inf\limits_{x\in A^2_H}I(x) =0 \label{E1}\\
 -\dot{u}(t) &=& \bar \partial
   L\big(t,u(t))\quad\forall t\in [0,T]\label{E2}\\
 -\frac{u(0)+u(T)}{2}  & \in& \bar \partial \ell (u(0)-u(T)).\label{E3}
\end{eqnarray}
\end{proposition}

\noindent{\bf Proof:}  Define for each $\la>0$, the $\lambda$-regularization $\ell^1_\la$ of the boundary Lagrangian $\ell$. By Lemma \ref{2-reg},  $\ell^1_\la$ is also anti-seldual on $H\times H$ and by Proposition \ref{lift}, the Lagrangian 
\begin{eqnarray*}
{\cal M}_\lambda (u,p)=\int_0^T L\big( t,u(t)+p_0(t),\dot u(t)\big)\, dt
  +\ell^1_\lambda \left( u(0)-u(T)+p_1,\frac{u(0)+u(T)}{2}\right)
\end{eqnarray*}
is anti-selfdual Lagrangian on $A_H^2\times (L_H^2\times H)$. It also satisfies  the hypothesis of Theorem \ref{one}. It follows that the infimum of the functional 
\[
I_\la(x)=\int_0^T L\big( t, x(t), \dot{x}(t)\big)\,
dt
  +\ell^1_\la \big(x(0)-x(T),
  \frac{x(0)+ x(T)}{2}\big)
\]
on $A^2_H$ is zero and is attained  at some  $x_\la \in A_H^2$ satisfying:
\begin{eqnarray}\label{  principle lambada }
\int_0^T  L\big( t, x_\la(t),\dot x_\la (t)\big)\, dt
  &+&\ell^1_\la \big( x_\la (0)-x_\la (T),\frac{x_\la (0)+x_\la (T)}{2}\big)
  =0\\
  -\dot x_\la(t)
  &\in& \bar \partial L\big( t,x_\la (t)) \\
 - \frac{x_\la (0)+x_\la (T)}{2}
  &\in& \bar \partial \ell^1_\la\big( x_\la (0)-x_\la (T)
\big).
\end{eqnarray}
We now show that $ (x_\la)_\la$ is bounded  in $A^2_H$.  Indeed,
since $\ell$ is bounded from below,  so is $\ell_\la$,  which together with
(\ref{  principle lambada })  imply that $\int_0^T L\big( t,
x_\la(t),\dot x_\la (t)\big)\, dt$ is bounded. It follows from
$(A'_1)$ and Lemma \ref{mild.estimate}  that $\{\dot x_\la (t)\}_\la$ is bounded in
$L^2_H$. It also follows from $(A'_2)$ that $\{ x_\la (t)\}_\la$ is
bounded in $L^2_H$,  hence, $x_\la$ is bounded in $A_H^2$ and thus, up
to a subsequence $ x_\la (t)\rightharpoonup u(t)$ in $A_H^2$,
$x_\la (0)\rightharpoonup u(0)$ and $x_\la (T)\rightharpoonup
u(T)$ in $H$.

From (\ref{ principle lambada }), we have that $ \ell^1_\la\big( x_\la
(0)-x_\la (T),\frac{x_\la (0)+x_\la (T)}{2}\big) $ is bounded from
above. Hence,  it follows from Lemma \ref{2-reg}   that
\begin{eqnarray*}
\ell \big( u(0)-u(T),\frac{u(0)+u(T)}{2}\big)
  \leq\liminf\limits_{\la \ra 0} \ell^1_\la
  \big( x_\la (0)-x_\la (T),\frac{x_\la (0)+x_\la (T)}{2}\big).
\end{eqnarray*}
By letting $\la\ra 0$ in (\ref{  principle lambada }),   we get
\begin{eqnarray*}
\int_0^T  L\big( t,u(t), -\dot {u}(t)\big)\, dt+
  \ell\big( u(0)-u(T),\frac{u(0)+u(T)}{2}\big)\leq 0.
\end{eqnarray*}
On the other hand, for every $x \in A^2_H$ we have
\begin{eqnarray*}
I(x)&=&\int_0^T  L\big( t,x(t),\dot x(t)\big)\, dt+
  \ell \big( x(0)-x(T),\frac{x(0)+x(T)}{2}\big)\\
  &=&\int_0^T \big\{L\big( t,x(t),\dot x(t)\big) +\langle x(t), \dot x (t)\rangle\big\} dt+
  \ell \big( x(0)-x(T),\frac{x(0)+x(T)}{2}\big)+\langle x(0)-x(T),\frac{x(0)+x(T)}{2}\rangle \\
 & \geq& 0
\end{eqnarray*}
which  means $I(u)=0$ and therefore 
$u(t)$ satisfies (\ref{E2}) and (\ref{E3}) as
well. $\square$

 \subsection{Parabolic semi-linear  without a diffusive term} 
We now consider the case where  $A$ is a purely skew-adjoint operator  and cannot therefore contribute to the coercivity of the problem.
 
 \begin{theorem} \label{no.diffusion} Let $(S_t)_{t\in \R}$ be a $C_0$-unitary group of operators associated to  a skew-adjoint operator $A$ on a Hilbert space $H$, and let   $\phi : [0, T]\times H\to \R\cup\{+\infty\}$ be a time-dependent convex,  Gateaux-differentiable function on $H$.
 Assume the following conditions:
\begin{description}
\item ($A_1)$ \quad For some $m,n >1 $ and $C_1, C_2>0$, we have for every $x\in L^2_H$,
 $$C_1\big(\| x\|_{L^2_H}^m-1\big) \leq \int_0^T \big\{\phi(t,x(t))+\phi^*(t, 0)\big\}\, dt\leq C_2 \big(1+\| x\|_{L^2_H}^n\big)$$  
\item ($A_2)$ \quad $\psi$ is a bounded below convex lower semi-continuous function on $H$ with $0 \in Dom (\psi).$
\end{description}
For any given $\omega\in\R$ and  $T >0$, consider the following functional on  $A^2_H,$
\[
I(x)=\int_0^Te^{-2\omega t}\left\{\phi(t, e^{\omega t}S_tx(t))+ \phi^*(t, -e^{\omega t}S_t\dot{x}(t))\right\}dt + \psi (x(0)-x(T))+\psi^*(-\frac{x(0)+ x(T)}{2}).
\]
Then, there exists a path  $ u 
\in A^p_H$  such that:
\begin{enumerate}
\item
$ I(u)=\inf\limits_{x\in A^p_H}I(x)=0$.
\item The path $v(t) := S_te^{\omega t}u(t)$ is a mild solution of the equation
\begin{eqnarray}
 \label{BC0}
\dot{ v}(t) + A v(t) +\omega v(t) &=&-\partial\phi( t, v(t))\quad \mbox{ for a.e. }t\in [0,T] \\
\frac{ v(0)+S_{-T}e^{-wT} v(T)}{2}
   &\in& - \partial \psi \left(  v(0)-S_{-T}e^{-wT} v(T)\right). \label{BC01}
\end{eqnarray}
Equation (\ref{BC0}) means that $v$ satisfies the following integral equation:
\begin{equation}\label{mild}
\hbox{\rm $ v(t)= S_tv(0)- \int_{0}^t S_{t-s}(\partial
\phi (s,v(s))-w v(s) \big ) \, ds$ for every $t\in [0, T]$.}
\end{equation}
\end{enumerate}
 \end{theorem}

\noindent {\bf Proof:}
Consider the anti-selfdual Lagrangians $M(t, x,p)=\phi (t, x)+\phi^*(t, -p)$ and $\ell (x,p)=\psi(
x)+\psi^*(-p)$, and  apply Proposition  \ref{prop.1}  to the Lagrangian
\begin{equation}
L(t, x, p)=e^{-2wt}  M\big(
t,S_te^{wt}x,S_te^{wt}p\big)
\end{equation}
which is anti-selfdual according to Proposition \ref{permanence}. 
 We  then obtain  ${u}(t)\in
A_H^2$ such that
\begin{eqnarray*}
\int_0^Te^{-2wt} \phi\big( t,S_te^{wt}u(t)\big) +\phi^*\big(-S_te^{wt}\dot{{u}}(t)\big)\, dt+\psi (u(0)-u(T)) +\psi^*(-\frac{u(0)+u(T)}{2}) =0,
\end{eqnarray*}
which gives
\begin{eqnarray*}
0 =&\int_0^T e^{-2wt}\left[\phi\big( t,S_te^{wt} u(t)\big) +\phi^*\big(-S_te^{wt}\dot{{u}}(t)\big) +\braket {S_te^{wt} u(t)}{S_te^{wt}\dot{{u}}(t)}\right]\, dt\\ &- \int_0^T \braket
{S_t u(t)}{S_t\dot{{u}}(t)} \,dt + \psi (u(0)-u(T)) +\psi^*(-\frac{u(0)+u(T)}{2}) \\
=&\int_0^T e^{-2wt}\left[\phi\big( t,S_te^{wt} u(t)\big) +\phi^*\big(-S_te^{wt}\dot{{u}}(t)\big) +\braket {S_te^{wt} u(t)}{S_te^{wt}\dot{{u}}(t)}\right]\, dt\\&- \int_0^T \braket
{u (t)}{\dot{{u}}(t)} \,dt + \psi (u(0)-u(T)) +\psi^*(-\frac{u(0)+u(T)}{2})  \\
=&\int_0^T e^{-2wt}\left[\phi\big( t,S_te^{wt} u(t)\big) +\phi^*\big(-S_te^{wt}\dot{{u}}(t)\big) +\braket {S_te^{wt} u(t)}{S_te^{wt}\dot{{u}}(t)}\right]\, dt\\& - \frac{1}{2} \|
 u(T) \|^2+ \frac{1}{2} \| u (0)\|^2+ \psi (u(0)-u(T)) +\psi^*(-\frac{u(0)+u(T)}{2}) \\
  =&\int_0^Te^{-2wt}\left[\phi\big( t,S_te^{wt}u(t)\big) +\phi^*\big(-S_te^{wt}\dot{{u}}(t)\big) +\braket
{S_te^{wt}u(t)}{S_te^{wt}\dot{{u}}(t)}\right]\, dt \\&+
\braket {u(0)-u(T)}{\frac{u(0)+u(T)}{2}}+ \psi (u(0)-u(T)) +\psi^*(-\frac{u(0)+u(T)}{2}) .
\end{eqnarray*}
Since  clearly 
$
 \phi\big( t,S_te^{wt} u(t)\big) +\phi^*\big(-S_te^{wt}\dot{{u}}(t)\big) +\braket {S_te^{wt} u(t)}{S_te^{wt}\dot{{u}}(t)}\geq 0
$ for every $t\in [0, T]$
and since 
$ \psi (u(0)-u(T)) +\psi^*(-\frac{u(0)+u(T)}{2}) +  \braket {u(0)-u(T)}{\frac{u(0)+u(T)}{2}} \geq 0
$,  we get equality from which we can conclude that 
\begin{eqnarray}
\hbox{$-S_t e^{wt}\dot {u}(t)= \partial \phi (t, S_te^{wt}u(t))$ for almost all $t\in [0, T]$ and $\frac{u(0)+u(T)}{2}\in -\partial \psi(u(0)-u(T))$.}
\end{eqnarray}
 In order to show  that $v(t):=S_t e^{wt}u(t)$ is a mild solution for (\ref{BC0}), we  set $x(t)=e^{wt}u(t)$ and write 
\[-S_t(\dot x(t)-w x(t))= \partial \phi (t, S_tx(t)),\] hence
$-(\dot x(t)+w x(t))= S_{-t}\partial \phi (t, v(t))$. 
By integrating between $0$ and $t$, we get  
\begin{eqnarray*}
x(t)=x(0)- \int_{0}^t \left\{S_{-s}\partial \phi (s,v(s))-w u(s) \right\} \, ds
\end{eqnarray*}
Substituting  $v(t)=S_t x(t)$ in the above equation gives
\begin{eqnarray*}
S_{-t}v(t)=v(0)- \int_{0}^t S_{-s}\big(\partial \phi (s,v(s))-w x(s)
\big ) \, ds,
\end{eqnarray*}
and consequently
\begin{eqnarray*}
 v(t)=S_tv(0)-S_t \int_{0}^t \big (S_{-s}(\partial \phi (s,v(s))-w
v(s) \big ) \, ds= S_tv(0)- \int_{0}^t S_{t-s}\big(\partial
\phi (s,v(s))-w v(s) \big ) \, ds
\end{eqnarray*}
which means that $v(t)$
 is a mild solution for (\ref{BC0}). 
 
 On the other hand, it is clear that the boundary condition $\frac{u(0)+u(T)}{2}\in -\partial \psi(u(0)-u(T))$ translates after the change of variables into 
  \begin{eqnarray*} \frac{ v(0)+e^{-wT}
S(-T)v(T)}{2}
   \in - \partial \psi \left( v(0)- e^{-wT} S(-T)v(T)\right)
\end{eqnarray*}
and we are done.\hfill $\square$
 \subsubsection*{Example 1: The complex Ginzburg-Landau equations in $\mathbb{R}^N$}
As an illustration, we consider  the  following evolution on $\mathbb{R}^N$
\begin{eqnarray}
\dot u (t)+i\triangle u+\partial \varphi (t, u(t))+w u(t)=0 \quad \mbox{ for } t\in [0,T].
\end{eqnarray}
(1)\, Under the  condition:
  \begin{eqnarray}\label{mild.boundedness1}
C_1 ( \int_0^T \|u(t)\|^2_2 \, dt-1) \leq \int_0^T \varphi (t, u(t))\, dt \leq C_2 ( \int_0^T \|u(t)\|_2^2 \, dt+1)
\end{eqnarray}
 where  $C_1,C_2>0$, Theorem \ref{no.diffusion}  yields  a solution of
\begin{eqnarray}\left\{ \begin{array}{lcl}
\dot u (t)+i\triangle u+\partial \varphi (t, u(t)) +\omega u(t)&=&0\quad \mbox{ for } t\in [0,T] \\
\hfill e^{-wT}e^{-iT\triangle }u(T) &=&u(0).\\
\end{array}\right.
\end{eqnarray}
 (2) If $w\geq 0$, then one can replace $\varphi$ with the convex function $\Phi (x)=\varphi (x)+\frac {w}{2}\|x\|^2$ to obtain  solutions such that 
 \begin{equation}
 \hbox{$u(0)= e^{-iT\triangle }u(T)$ or $ u(0)=- e^{-iT\triangle }u(T)$.}
 \end{equation}
(3)\, One can also drop the coercivity condition  (the lower bound) on $\varphi (t, u(t))$ in (\ref{mild.boundedness1}) and still get  periodic-type solutions. Indeed, by applying our result to  the now coercive convex functional $\Psi (t, u(t)):=\varphi (t, u(t))+\frac{\epsilon}{2}\| u(t)\|_H^2$, and $w-\epsilon$,   to obtain a solution such that
 \begin{eqnarray}
 \hfill e^{(-w+\epsilon) T}e^{-iT\triangle }u(T)&=&u(0).
\end{eqnarray}

\subsubsection*{Example 2:  Almost Periodic solutions for linear Schr\"odinger equations:}
Consider now the following linear Schrodinger equation
\begin{equation}
i\frac{\partial u}{\partial t}=\Delta u -V(x)u.
\end{equation}
Assuming that the space $ \{ u \in H^{2,2}(\mathbb{R}^N) : \int_{\R^N} |V(x)|u^2 \, dx <\infty\}$ is
dense in $H:=L^{2}(\mathbb{R}^N)$, we get that the operator $Au:=i\triangle u-iV(x)u$ is skew adjoint on $H$. In order to introduce some coercivity, and to avoid the trivial solution, we can consider for any $\epsilon, \delta >0$ and $0 \neq f \in H$,  the convex function $\varphi_{\epsilon}  (u) := \frac {\epsilon}{2}
\|u\|_H^2+\delta \braket {f}{u}$.

By applying  Theorem \ref{no.diffusion} to $A$, $\phi_\epsilon$, and $\omega =\epsilon$, we get a non trivial
 solution $u\in A^2_H$ for the equation
\begin{eqnarray}\left\{ \begin{array}{lcl}
i\frac{\partial u}{\partial t}= \triangle u -V(x)u + \delta f,\\
u (0)=e^{-\epsilon T}e^{iT(-\triangle+V(x)) }u(T).\\
\end{array}\right.
\end{eqnarray}

\subsubsection*{Example 3:  Coupled flows and wave-type equations}
Let  $A:
D(A) \subseteq H \rightarrow H$ is a linear operator with a dense
domain in $H.$ Suppose $D(A)=D(A^*)$, and define the following operator ${\cal A} $ on the product space $H\times H$ as follows:
\begin{eqnarray*}\left\{ \begin{array}{lcl}
{\cal A}: D({\cal A} ) \subseteq H \times H \rightarrow H \times H,\\
{\cal A}  (x,y):= (Ay, -A^* x)\\
\end{array}\right.
\end{eqnarray*}
It is easily seen that  ${\cal A}: D({\cal A}) \subseteq H \times H
\rightarrow H \times H$ is a skew-adjoint operator, and hence by
 virtue of Stone's Theorem,  ${\cal A} $ is the generator of a  $C_0$ unitary group $\{S_t\}$ on $H \times H.$
 Here is another application of Theorem \ref{no.diffusion}.
 \begin{theorem} Let  $ \varphi (t, \cdot)$  and $\psi$ be proper convex lower semi continuous functionals  on $H \times H.$   Asume the following conditions:
\begin{description}
\item ($A''_1)$ \quad For some $ m,n >1 $ and $C_1, C_2>0$, we have
 \[
 \hbox{$C_1\big(\| x\|_{L^2_{H \times H}}^m-1\big) \leq  \int_0^T \varphi (t,x(t))\, dt\leq C_2 \big(1+\| x\|_{L^2_{H \times H}}^n\big)$ for
  every $x\in L^2_{H \times H}$.}
  \] 
\item ($A''_2)$ \quad $\psi :H\times H\ra \R \cup\{ +\infty\}$  is  bounded  below and $0 \in Dom (\psi).$
\end{description}
Then  there exists 
a mild solution $(u(t), v(t))\in A^2_{H \times H}  $
for the following system,
\begin{eqnarray*}\left\{ \begin{array}{lcl}
 -\dot {u}(t) +A v(t)+w u(t)=\partial _1\varphi (t,u(t),v(t)),\\
-\dot {v}(t) -A^*u(t)+w v(t)=\partial _2\varphi (t,u(t),v(t))\\
\end{array}\right.
\end{eqnarray*}
with a  boundary condition of the form  (\ref{BC01}).
\end{theorem}
 \subsection{Parabolic semi-linear equation with a diffusive term}  
The existence of periodic solutions follows from a more general result of Lions
(See [12; Proposition III.5.1]).
 Our approach is quite different and relies on last section's selfdual variational principle which will now yield  true periodic solutions provided the strong coercivity conditions of Lions are satisfied.   
 
For  given $0<T<\infty$,  $1< p< \infty$, and a Hilbert space $H$ such that  $X\subseteq H \subseteq X^*$ is an evolution triple, we  consider  the space
\begin{eqnarray*}
{\cal X}_{p,q}=\{ u: u  \in L^p(0,T:X),   \dot u \in
L^q(0,T:X^*) \}
\end{eqnarray*}
equipped with the norm
$\|u\|_{{\cal X}_{p,q}}= \|u\|_{L^p(0,T:X)}+\| \dot u\|_{L^q(0,T:X^*)}, $
which leads to a continuous injection
${\cal X}_{p,q} \subseteq C( 0,T:H)$.
We shall  prove the following.

   \begin{theorem}\label{diffusion} Let $X\subset H \subset X^*$ be an evolution triple, and consider a time-dependent anti-selfdual Lagrangian $L(t,x,p)$ on
$[0,T] \times X \times X^*$ and an anti-selfdual Lagrangian $\ell$ on
$H\times H$ such that the following conditions are satisfied:
\begin{description}
\item ($B'_1)$ \quad For some $p\geq 2, m,n >1 $ and $C_1, C_2>0$, we have
 \[
 \hbox{$C_1\big(\| x\|_{L^p_X}^m-1\big) \leq  \int_0^T L(t,x(t),0)\, dt\leq C_2 \big(1+\| x\|_{L^p_X}^n\big)$  for every $x\in L^p_X$.}
 \]
\item ($B'_2)$ \quad $\ell$ is bounded from below. 
\end{description}
The  following functional   
\[
I(x)=\int_0^T   {L}\big( t, x(t),\dot{ x}(t)\big)\, dt
  +\ell\big( x(0)- x(T),
  \frac{ x(0)+ x(T)}{2}\big)
\]
then attains its minimum on ${\cal X}_{p,q}$ at a path  $ u
\in {\cal X}_{p,q}$   such that
 \begin{equation}\label{principle 1}
 I(u)=\inf \{I(x); x\in {\cal X}_{p,q}\}=0
 \end{equation}
  \begin{equation} \label{principle 2}
-\dot{{u}}(t) \in  \bar \partial { L} \big(t,u(t)\big)\quad\forall t\in [0,T]
   \end{equation}
 \begin{equation}\label{principle 3}
-\frac{u(0)+u(T)}{2}
  \in \bar  \partial \ell \big(u(0)-u(T)\big).
\end{equation}
\end{theorem}
 {\bf Proof:} Use Lemma \ref{mild.estimate} to lift the Lagrangian $L$ to a time dependent  ASD Lagrangian on $[0, T]\times H\times H$ via the formula 
 \begin{eqnarray*}
M(t, u,p):=
\left\{ \begin{array}{lcl}
L(t, u,p), \quad \quad \quad u\in X,\\
+\infty \quad \quad \quad u\in H \setminus X.
\end{array}\right.
\end{eqnarray*}
 We start by assuming that 
$\ell (a,b)\rightarrow \infty $ as $\|b\|\rightarrow \infty$.  Consider 
for  $\la
>0$, the $\la-$regularization of $M$, namely
\begin{eqnarray}\label{regularization.1}
L^1_\la (t,x,p):=\inf\left\{ M(t,z,p)+\frac{\|
x-z\|^2}{2\la}+\frac{\la}{2}\| p\|^2; \, z\in H\right\}.
\end{eqnarray}
It is easy to check that $L_\la$ satisfies the conditions $(A'_1)$ and
$(A'_2)$ of Proposition \ref{no.diffusion}. It  follows  that there exists 
 a path  $x_\la (t)\in A_H^2$  that
\begin{eqnarray} \label{regularization.2}
 \int_0^T  L^1_{\la}\big( t, x_\la(t),\dot x_\la (t)\big)\, dt
  +\ell \big( x_\la (0)-x_\la (T),\frac{x_\la (0)+x_\la (T)}{2}\big) =0.
 \end{eqnarray}
We now  show that $(x_\la)_\la$ is bounded in an appropriate
function space. Indeed, since $L$ is convex and lower semi-continuous, there
exists $i_\la (x_\la)$ such that the infimum in (\ref{regularization.1}) is attained at $i_\la (x_\la)\in X,$ i.e.
\begin{eqnarray}\label{regularization.3}
L_\la (t, x_\la (t), \dot {x}_\la (t) )=
L(t,i_\la (x_\la), \dot {x}_\la (t) )
+\frac{\| x_\la (t)- i_\la (x_\la)\|^2}{2\la}+ \frac{\la}{2}\|
 \dot {x}_\la (t)\|^2.
\end{eqnarray}
Plug (\ref{regularization.3}) in equality  (\ref{regularization.2})   to get
\begin{eqnarray}\label{regularization.4}
\int_0^T\big( L(t,i_\la (x_\la), \dot {x}_\la (t) )+\frac{\| x_\la (t)- i_\la (x_\la)\|^2}{2\la}
  +\frac{\la}{2}\|  \dot {x}_\la (t)\|^2  \big) dt  
  +\ell\big( x_\la (0)-x_\la (T),\frac{x_\la (0)+x_\la
  (T)}{2}\big)=0.
 \end{eqnarray}
By the coercivity assumptions in $(B'_1)$, we obtain that $(i_\la
(x_\la))_\la$ is bounded in $L^p(0,T; X)$ and $(x_\la)_\la$ is
bounded in $ L^2(0,T; H)$. According to Lemma \ref{mild.estimate},  Condition ($B'_1$) yields that 
$ \int_0^T L(t,x(t),p(t))\, dt$ is coercive in $p(t)$ on $L^q(0,T; X^*)$, and therefore 
it follows from  (\ref{regularization.4})  that
  $ (\dot x_\la)_\la$ is bounded
in $L^q(0,T; X^*)$. Also, since $L$ and $\ell$ are bounded from below,
it follows again from (\ref{regularization.4}) that
$\int_0^T\| x_\la (t)- i_\la (x_\la)\|^2 \, dt  \leq 2\la C$ 
for a constant $C>0$. Since now $x_{\la} (T)-x_{\la}(0)=\int_0^T \dot
x_\la \, dt$,  therefore $x_{\la} (T)-x_{\la}(0)$ is bounded in
$X^*$. Also, since we have assumed that $\ell (a,b)\rightarrow \infty $ as
$\|b\|\rightarrow \infty,$ it follows that  $x_{\la} (0)+x_{\la}(T)$
is also bounded in $H$ and  consequently in $X^*$.  Therefore there
exists $u  \in L^2_H$ with  $\dot { u} \in L^q(0,T; X^*)$
and $u(0), u(T) \in X^*$  such that
\begin{eqnarray*}
i_\la (x_\la) \rightharpoonup  u \quad \text{  in  } L^p(0,T; X),\\
  \dot x_\la \rightharpoonup  \dot
{ u}\quad \text{  in  }  L^q(0,T; X^*),\\
{x_\la} \rightharpoonup  u \quad \text{  in  }  L^2(0,T; H)\\
x_\la(0)\rightharpoonup u(0), \quad x_\la(T)\rightharpoonup
u(T) \quad \text{  in  }  X^*.
\end{eqnarray*}
By letting $\lambda$ go to zero in (\ref{regularization.4}), we
obtain from the above that
\begin{equation}\label{negative}
 \ell \big( u (0)-u (T),\frac{u (0)+u (T)}{2}\big)
  + \int_0^T  L \big( t,   {u} (t), \dot {u}(t) \big) \, dt
    \leq 0.
\end{equation}
It follows from $(B'_1)$, Lemma \ref{mild.estimate} and (\ref{negative}) that
$u  \in {\cal X}_{p,q}$  and consequently, $u(0), u(T)
\in H$. \\  
Now we show that
 one can actually do without  the coercivity condition on $\ell.$
 Indeed, by using  the $\la-$regularization
$\ell^1_{\la}$ of $\ell$, we get the required coercivity condition on the
second variable of $\ell_{\la}$ and we obtain from the above
  that there exists  $ x_{\la} \in{\cal X}_{p,q}$  such that
\begin{eqnarray}\label{regularization.5}
\int_0^T  L\big( t, x_{\la}(t),\dot{x}_{\la}(t)\big)\, dt
  +\ell_{\la}\left( x_{\la}(0)- x_{\la}(T),
  \frac{ x_{\la}(0)+ x_{\la}(T)}{2}\right) \leq 0.
\end{eqnarray}
It follows from $(B'_1)$ and the boundedness of $\ell^1_{\la}$ from below, 
that $ (x_{\la})_\la$ is bounded in $ L^p(0,T; X)$, and $(\dot
{x}_{\la})_\la$ is bounded in $ L^q(0,T; X^*)$ again  by virtue of Lemma \ref{mild.estimate}. Hence, $(x_{\la})_\la$ is bounded in ${\cal X}_{p,q}$ and therefore
$(x_{\la}(0))_\la$ and $ (x_{\la}(T))_\la$ are bounded in $H$. We therefore get, up to a subsequence, that 
\begin{eqnarray*}
 x_\la \rightharpoonup u \quad \text{  in  } L^p(0,T; X),\\
  \dot x_\la \rightharpoonup  \dot
{ u}\quad \text{  in  }  L^q(0,T; X^*),\\
x_\la(0)\rightharpoonup u(0) \quad \text{  in  } H,\\
x_\la(T)\rightharpoonup u(T) \quad \text{  in  } H.
\end{eqnarray*}
By letting $\lambda$ go to zero in (\ref{regularization.5}), it
follows from  the above   that
\[
  \int_0^T   L \big( t,    {u} (t), \dot {u}(t) \big) \,
  dt+\ell\big( u (0)-u (T),\frac{u (0)+u (T)}{2}\big)
    \leq 0.\]
So $I(u)=0$ and $u$ is a solution of (\ref{principle 2}) and (\ref{principle 3}).

 \begin{corollary} \label{theorem.2}  Let $X\subseteq H \subseteq X^*$ be  an evolution triple,  Let $A:X\to X^*$  a bounded positive operator on $X$ and let $\phi : [0, T]\times X\to \R\cup\{+\infty\}$ be a time-dependent convex, lower semi-continuous and proper function on $X$.  Consider the convex function $\Phi (x)=\phi (x) +\frac{1}{2} \langle Ax, x \rangle$ as well as the anti-symmetric part $A^a:=\frac{1}{2}(A-A^*)$ of $A$.
Assume the following conditions hold:
\begin{description}
\item ($B_1)$ \quad For some $p\geq 2, m,n >1 $ and $C_1, C_2>0$, we have for every $x\in L^p_X$,
 $$C_1\big(\| x\|_{L^p_X}^m-1\big) \leq \int_0^T \left\{\Phi(t,x(t))+\Phi^*(t, -A^ax(t)\right\} dt\leq C_2 \big(1+\| x\|_{L^p_X}^n\big)$$
\item ($B_2)$ \quad $\psi$ is a bounded below convex lower semi-continuous function on $H$ with $0 \in Dom (\psi).$
\end{description}
For any  $T >0$ and $\omega\geq 0$,  consider the following functional on  ${\cal X}_{p,q}$
\[
I(x)=\int_0^Te^{-2\omega t}\left\{\Phi(t, e^{\omega t}x(t))+ \Phi^*(t, -e^{\omega t}(A^ax(t) +\dot{x}(t)))\right\}dt + \psi (x(0)-x(T))+\psi^*(-\frac{ x(0)+ x(T)}{2}).
\]
Then, there exists a path  $ u
\in L^p(0,T:X)$  with  $ \dot {{u}} \in L^q(0,T:X^*)$ such that:
\begin{enumerate}
\item
$ I(u)=\inf\limits_{x\in {\cal X}_p}I(x)=0$.
\item If $ v(t)$ is defined by $v(t) := e^{\omega t}u(t)$ then it satisfies
\begin{eqnarray}
 \label{eqn:1000}
\dot{ v}(t) +Av(t)+\omega v(t) &\in&-\partial\phi( t, v(t))\quad \mbox{ for a.e. }t\in [0,T] \\
\frac{ v(0)+e^{-wT} v(T)}{2}
   &\in& - \partial \psi \left( v(0)-e^{-wT} v(T)\right). \label{BC02}
\end{eqnarray}
\end{enumerate}
 \end{corollary}

\noindent{\bf Proof:}  It suffices to apply  Theorem \ref{diffusion}   to
the  anti-selfdual Lagrangian 
\[
L(t, x,p)=e^{-2\omega t}\left\{ \Phi (t, e^{\omega t}x)+\Phi^*(t, -e^{\omega t}A^ax
-e^{\omega t}p)\right\}
\]
 associated to a convex lower semi-continuous function $\Phi$, a
skew-adjoint operator $A^a$ and a scalar $\omega$.   $\square$

\subsubsection*{Example 4: Complex Ginzburg-Landau evolution with diffusion}

Consider a complex Ginzburg-Landau equations of the following type.

\begin{eqnarray}\label{GL.19}\left\{ \begin{array}{lcl}
\frac{\partial u}{\partial t} -(\kappa+i)\triangle u+\partial \Psi (t,u)+wu&=&0, \quad \quad (t,x) \in (0,T) \times \Omega,\\
\hfill u(t,x)&=&0 \quad \quad \quad \quad x  \in \partial  \Omega,\\
\hfill e^{-wT}u(T)&=&u(0),\\
\end{array}\right.
\end{eqnarray}
where $\kappa > 0$, $\omega \leq 0$, $\Omega$ is a bounded domain in $\mathbb{R}^N$ and $\Psi$ is a time-dependent convex lower semi-continuous function. An immediate corollary of Theorem \ref{diffusion} is the following.

\begin{corollary}  \label{GL.20} Let $X:= H_0^{1}(\Omega)$, $H:=L^2(\Omega)$ and $X^*=H^{-1} (\Omega).$ If  for some $C>0$, we have
\begin{eqnarray*}
\hbox{$-C\leq \int_0^T \Psi (t, u(t))\, dt \leq C ( \int_0^T \|u(t)\|_{H^1_0}^2 \, dt+1)$ for every $u \in L^2_X,$}
\end{eqnarray*}
  then there exists a solution $u \in  {\cal X}_{2,2}
$ for (\ref{GL.19}).
\end{corollary}
 
\noindent{\bf Proof}: Set $\varphi (t,u):= \frac
{k}{2}    \int_\Omega |\nabla u|^2 dx+\Psi (t, u(t))$, $A=-(1+i)\Delta$, $A^a=-i\Delta$  and
note that since  \begin{equation}\label{bound.1} c_1
(\|u\|^2_{L^2_X}-1)\leq \int_0^T \varphi (t,u) \, dt \leq c_2
(\|u\|^2_{L^2_X}+1)
\end{equation}
for some $c_1, c_2>0,$ we therefore have
\begin{eqnarray*}
c'_1 (\|v\|^2_{L^2_{X^*}}-1)\leq \int_0^T\varphi^* (t,v) \, dt \leq c'_2 (\|v\|^2_{L^2_{X^*}}+1)
\end{eqnarray*}
for some $c'_1, c'_2>0$, and hence
\begin{eqnarray*}
c'_1 (\int_0^T \int_{\Omega}|\nabla (-\triangle)^{-1}v|^2\, dx \,
dt-1)  \leq \int_0^T\varphi^* (t,v) \, dt \leq c'_2  (\int_0^T
\int_{\Omega}|\nabla (-\triangle)^{-1}v|^2\, dx \, dt+1).
\end{eqnarray*}
from which we obtain
\begin{eqnarray*}
c'_1 (\int_0^T \int_{\Omega}|\nabla u|^2\, dx \, dt-1)  \leq \int_0^T \varphi^* (t,i\triangle u) \, dt \leq c'_2  (\int_0^T \int_{\Omega}|\nabla u|^2\, dx \, dt+1)
\end{eqnarray*}
which,  once coupled with (\ref{bound.1}), yields the required boundedness in  $(B'_1)$. \hfill $\square$\\

We now show how one  can sometimes combine the two ways to define an ASD Lagrangian that deals with a superposition of an unbounded skew adjoint operators with another bounded positive operator. Note the impact on the boundeness condition  $(B_1)$ above.

\begin{corollary} \label{ntheorem.2}  Let $X\subseteq H \subseteq X^*$ be  an evolution triple in such a way that the duality map $D:X \to X^*$ is linear and symmetric. Let $A_1:X\to X^*$ be  a bounded positive operator on $X$ and let $A_2: D(A) \subseteq X\to X^*$ be a --possibly unbounded--  skew adjoint operator. Let $\phi : [0, T]\times X\to \R\cup\{+\infty\}$ be a time-dependent convex, lower semi-continuous and proper function on $X$, and consider the convex function $\Phi (x)=\phi (x) +\frac{1}{2} \langle A_1x, x \rangle$ as well as the anti-symmetric part $A_1^a:=\frac{1}{2}(A_1-A_1^*)$ of $A_1$. Let $\bar S_t: X^*\to X^*$ be the unitary group generated by $A_2$ and $S_t=D \bar S_t D^{-1}.$
Assume the following conditions:
\begin{description}
\item[ ($D_1$)]  For some $p\geq 2, m,n >1 $ and $C_1, C_2>0$, we have for every $x\in L^p_X$,
 $$C_1\big(\| x\|_{L^p_X}^m-1\big) \leq \int_0^T \left\{\Phi(t,  S_t x(t))+\Phi^*(t, -A_1^a  S_t x(t)\right\} dt\leq C_2 \big(1+\| x\|_{L^p_X}^n\big)$$ 
\item[ ($D_2$)]   $\psi$ is a bounded below convex lower semi-continuous function on $H$ with $0 \in Dom (\psi).$
\end{description}
For any  $T >0$ and $\omega\in \R$,  consider the following functional on  ${\cal X}_{p,q}$
\begin{eqnarray}
I(x)&=&\int_0^Te^{-2\omega t}\left\{\Phi(t, e^{\omega t}  S_t x(t))+ \Phi^*(t, -e^{\omega t}(A^a S_t x(t) + \bar S_t \dot{x}(t)))\right\}dt\\
&& + \psi (x(0)-x(T))+\psi^*(-\frac{ x(0)+ x(T)}{2}).
\end{eqnarray}
Then, there exists a path  $ u
\in L^p(0,T:X)$  with  $ \dot {{u}} \in L^q(0,T:X^*)$ such that:
\begin{enumerate}
\item
$ I(u)=\inf\limits_{x\in {\cal X}_p}I(x)=0$.
\item Moreover, if  $ \bar S_t=  S_t$ on $X$, then $v(t) := e^{\omega t} \bar S_t u(t)$  satisfies
\begin{eqnarray}
 \label{neqn:1000}
\dot{ v}(t) +A_1v(t) +A_2v(t) +\omega v(t) &\in&-\partial\phi( t, v(t))\quad \mbox{ for a.e. }t\in [0,T] \\
\frac{ v(0)+S_{(-T)}e^{-wT} v(T)}{2}
   &\in& - \partial \psi \left( v(0)-S_{(-T)}e^{-wT} v(T)\right). 
\end{eqnarray}
\end{enumerate}
 \end{corollary}
\noindent{\bf Proof:}  It suffices to apply  Theorem \ref{diffusion}   to
the  anti-selfdual Lagrangian 
\[
L_S(t, x,p)=e^{-2\omega t}\left\{ \Phi (t, e^{\omega t}S_t x)+\Phi^*(t, -e^{\omega t}A^a S_t x
-e^{\omega t}\bar S_t p)\right\}
\]
 which is anti-selfdual in view of the remark of section 2. \hfill  $\square$
\subsubsection*{Example 5: The complex Ginzburg-Landau equations with advection in a bounded domain}
We consider  the  following evolution on bounded domain  $\Omega,$
\begin{eqnarray*}
\dot u (t)-i\triangle u+a.\nabla u(t)+\partial \varphi (t, u(t))+ w u(t)=0 \quad \mbox{ for } t\in [0,T].
\end{eqnarray*}
 Under the  condition that $a$ is a constant vector and 
  \begin{eqnarray}\label{nmild.boundedness1}
C_1 ( \int_0^T \|u(t)\|^2_2 \, dt-1) \leq \int_0^T \varphi (t, u(t))\, dt \leq C_2 ( \int_0^T \|u(t)\|_2^2 \, dt+1)
\end{eqnarray}
 where  $C_1,C_2>0$, Corollary  \ref{ntheorem.2}  yields  a solution of
\begin{eqnarray}\left\{ \begin{array}{lcl}
\dot u (t)-i\triangle u+a.\nabla u(t)+\partial \varphi (t, u(t)) +\omega u(t)&=&0\quad \mbox{ for } t\in [0,T] \\
\hfill e^{-wT}e^{-iT\triangle }u(T) &=&u(0).\\
\end{array}\right.
\end{eqnarray}
\noindent{\bf Proof}: Set $A_1u=a \cdot \nabla u$, $A_2=-i\Delta$  and $H=L^2(\Omega)$ in Corollary \ref{ntheorem.2}. Define the Banach space 
$X_1=\{u \in H; A_1u \in H  \}$
equipped with the norm $\|u\|_X=\big (\|u\|_H^2+\|A_1u\|_H^2\big)^{\frac{1}{2}}$.
Therefore $X^*=\{(I+A_1^*A_1)u ; u \in X \}$ 
and the norm in $X^*$ is
$\|f\|_{X^*}=\|(I+A_1^*A_1)^{-1}f\|_X$. 
Note that $D=I+A_1^*A_1$ is the  duality map between $X$
and $X^*$, since
$\langle u,Du\rangle=\langle u,(I+A_1^*A_1)u\rangle=\|u\|_X^2$
and
$\|D u\|_{X^*}=\|D^{-1}Du\|_X=\|u\|_X$.

\section{Hamiltonian systems with general boundary conditions}
In this section we consider the system
\begin{eqnarray}\label{1000}
J\dot u(t)+J {\cal A}u(t)&= &\bar \partial L\big( t,u(t)\big),
\end{eqnarray}
where $L$ is a time dependent anti-selfdual Lagrangian on $[0,T] \times X \times X$, where $X:=H \times H$ for some --possibly  infinite-- dimensional  Hilbert space $H$, ${\cal A}(p,q)= (Ap,-Aq)$
where $A:D(A)\subseteq H\rightarrow H$ is a  self-adjoint  operator, and  $J$ is the symplectic operator $J(p,q)=(-q,p)$.

We assume that  $\langle Au, u\rangle \geq c_0\|u\|^2_H$ on $D(A)$ for some $c_0>0$, and that $A^{-1}$ is compact. We shall denote by  $\tilde A$ the operator $(A, A)$ on the product space $X=H\times H$,   and consider the Hilbert space $Y \subseteq X$ which is the completion of $D(\tilde A)$ for the norm induced by the inner product
$
\langle u,v\rangle_{Y}:=\langle u, \tilde A v \rangle_X. 
$
   The path space 
 \[
 W= \big \{ u \in L^2_X[0, T]; \, \dot u \,  \& \,  \tilde A u \in L^2_X[0, T] \big \}
 \]
 is also a Hilbert space once equipped
 with the norm
$\|u\|_W= \big ( \|\tilde A_1 u\|^2_{L^2_X}+\|\dot u\|_{L^2_X}\big)^{\frac{1}{2}}.$
The embedding $W
\to  C([0,T];X)$ is then continuous, i.e.,
\begin{eqnarray}
\|u\|_{C([0,T];X)}\leq c\|u\|_W,
\end{eqnarray}
for some constant $c>0$, while the injection $W \to L^2([0,T];X)$ is compact.  We shall also consider  (\ref{1000}) with a  general boundary condition such as
\begin{equation}
 R\frac {v(T)+v(0)}{2}\in   \partial
\psi \big(v(T)-v(0)\big), 
\end{equation}
where $\psi$ is a convex lsc function on $X$ and $R$ is the  automorphism $R (p,q)=(p,-q)$ on $X$.\\
Here is our main variational principle for Hamiltonian systems.

\begin{theorem}\label{thm.new.var}Let $L:[0,T]\times X \times X \ra \R \cup \{+\infty \}$ be a time dependent anti-selfdual Lagrangian, and let $\ell:X \times X \ra \R \cup \{+\infty \}$ be a convex Lagrangian that satisfies
\begin{equation}
\hbox{ $\ell(x,p) \geq -\langle x, \tilde A p \rangle$ for all $x \in X$ and $p \in Y$.}
\end{equation}
Assume the following conditions:
\begin{description}
\item ($C_1$)\quad There exists $0< \beta  <\frac {1}{8c\sqrt{T}}$ and  $   \gamma , \al  \in L^2 (0,T;\mathbb{R_{+}})$
such that
\[
\hbox{$- \al (t)\leq L(t,u,0)\leq\frac{\beta}{2} \|u\|^2+\gamma (t)$ for every $u\in X$ and a.e. $t\in [0,T]$.}
\]
\item ($C_2$)\quad 
$0\in {\rm Dom}(\ell)$, $\ell$ is bounded below on $X$,  and its  restriction to $Y \times Y$ is lower semi-continuous for the $Y\times Y$ topology.   
\end{description}
 Then the infimum of the functional 
\begin{eqnarray} \label{new.var}
I(u)&=& \int_0^T \left\{ L (t,u(t),-J\dot u (t)-J{\cal A}u(t) )-\langle J \dot u(t)+J{\cal A}u(t),u(t)\rangle \right\} dt \nonumber\\
 &&+ \frac{1}{\beta} \ell \big( u(T)-u(0), R\frac {u(T)+u(0)}{2}\big) +( \frac{1}{\beta} \ell )^*\big( R\frac {u(T)+u(0)}{2}, u(T)-u(0)\big)\nonumber\\
&& -\langle u(T)-u(0),R( u(0)+u(T))\rangle
\end{eqnarray}
is  equal to zero and is attained at some $u \in W $ and  $u$ is  a solution of
\begin{eqnarray}\label{H-equation1}\left\{ \begin{array}{lcl}
\hfill J\dot v(t)+J{\cal A}v(t)&= &\bar \partial L\big( t,v(t)\big)\\
\hfill R\frac {v(T)+v(0)}{2}&= & \bar \partial \frac { \ell}{\beta} \big(v(T)-v(0)\big).
\end{array}\right.
\end{eqnarray}
\end{theorem}
We start by establishing the following proposition which assumes   a stronger condition on
the main Lagrangian $L$ and the boundary Lagrangian $\ell$.
\begin{proposition} \label{Relax.2} Let $L:[0,T]\times X \times X \ra \R \cup \{+\infty \}$,  $\ell: X \times X\ra \R \cup \{+\infty \}$ be as in Theorem \ref{thm.new.var}, and assume the following conditions:
\begin{description}
\item ($C'_1$)\quad  There exists $\la >0$  and $0< \beta  <\frac {1}{8c\sqrt{T}},$ and  $\gamma, \al  \in L^2 (0,T;\mathbb{R_{+}})$
such that
\[
\hbox{$- \al (t)\leq L(t,u,p)\leq\frac{\beta}{2} \|u\|^2+ +\la\|p\|^2+\gamma (t)$ for every $ (u,p) \in X \times X$ and a.e. $t\in [0,T]$.}
\]
\item ($C'_2$)\quad  There exist positive constants   $ {\al}_{1}, \beta_{1}, \gamma_{1} \in \mathbb{R}$ such that,
for every $ (u,v) \in Y \times Y$  one has
$$- \al_1\leq \ell(u,v) \leq\frac{\beta_1}{2}(\|u\|_Y^2+\|v\|_Y^2)+ \gamma_1.$$
\end{description}
 Then the functional 
\begin{eqnarray} \label{new.var}
I(u)&=& \int_0^T \left\{ L (t,u(t),-J\dot u (t)-J{\cal A}u(t) )-\langle J \dot u(t)+J{\cal A}u(t),u(t)\rangle \right\} dt \nonumber\\
 &&+ \frac{1}{\beta} \ell \big( u(T)-u(0), R\frac {u(T)+u(0)}{2}\big) +( \frac{1}{\beta} \ell )^*\big( R\frac {u(T)+u(0)}{2}, u(T)-u(0)\big)\nonumber\\
&& -\langle u(T)-u(0),R( u(0)+u(T))\rangle
\end{eqnarray}
is selfdual on $W$ and its corresponding anti-symmetric Hamiltonian on $W \times W$  is
\begin{eqnarray*}
 M(u, v)&=&\int_0^T\left\{\langle J\dot v (t)+J{\cal A}v(t) ,u(t)\rangle+\tilde H_{L}(t,-J\dot u (t)-J{\cal A}u(t),-J\dot v-J{\cal A}v(t)
(t))\right\}dt\\
&&-\int_0^T\langle J\dot u(t)+J{\cal A}u(t),u(t)\rangle\, dt 
 -\left\langle u(T)-u(0),R(u(T)+u(0))\right\rangle \\ &&+\left\langle u(T)-u(0),R\frac{v(T)+v(0)}{2}\right\rangle
+\left\langle v(T)-v(0),R\frac{u(T)+u(0)}{2}\right\rangle \\ 
&&   +\frac{1}{\beta}\ell \big( u(T)-u(0), R\frac
{u(T)+u(0)}{2}\big)-\frac{1}{\beta}\ell \big( v(T)-v(0), R\frac {v(T)+v(0)}{2}\big).
\end{eqnarray*}
Also, the infimum of $I$
 is  equal to zero and is attained at some $u \in W $ and  $u(t)$ is  a solution of
\begin{eqnarray}\label{H-equation01}\left\{ \begin{array}{lcl}
\hfill J\dot v(t)+J{\cal A}v(t)&= &\bar \partial L\big( t,v(t)\big)\\
\hfill R\frac {v(T)+v(0)}{2}&= & \bar \partial \frac {\ell}{\beta} \big(v(T)-v(0)\big).
\end{array}\right.
\end{eqnarray}
\end{proposition}
The proof requires a few preliminary  lemmas.  We first  establish the self duality of the functional $I$.

\begin{lemma}\label{H-positive1} With the above notation we have
\begin{enumerate}
\item For   every $u\in W$, we have  $I(u)\geq 0$. 
\item  $M$ is an anti-symmetric Hamiltonian on $W\times W$.
\item For every  $u\in W$, we have 
$
I(u)= \sup\limits_{v\in W}{M}(u, v).
$
\end{enumerate}
 \end{lemma}
\paragraph{Proof:} 1)\, 
Since $L$ is anti-selfdual Lagrangian we have or any $u\in W$,
\begin{eqnarray*}
L(t,u(t),-J\dot u(t)-J{\cal A}u(t))-\langle J\dot
u(t)+J{\cal A}u(t),u(t)\rangle \geq 0 \quad  \quad \text{for }t \in [0,T].
\end{eqnarray*}
Also it follows from the definition  of   Legendre-Fenchel duality that 
\begin{eqnarray*}
\frac {1}{\beta}\ell \big( u(T)-u(0), R\frac {u(T)+u(0)}{2}\big) +(\frac {1}{\beta}\ell)^*\big( R\frac {u(T)+u(0)}{2}, u(T)-u(0)\big) -\langle u(T)-u(0),R( u(0)+u(T))\rangle \geq 0.
\end{eqnarray*}
from which we obtain  $I(u)\geq 0$. \\
2)\,  The fact that $M$ is an anti-symmetric Hamiltonian on $W\times W$ is straightforward. Indeed, the weak lower semi-continuity of $u\to M(u, v)$ for any $v\in W$ follows from the fact that  the embedding $W\subseteq L^2_X$ is compact  and $W\subseteq C(0,T;X)$ is continuous.  It follows that if  $u \in W$ and $\{u_n\}$ is  a bounded sequence in $W$ such that $u_n \rightharpoonup u$ weakly in $W$, then 
\begin{eqnarray*}
\lim_{n\rightarrow \infty}\int_0^T \langle J\dot u_n+J{\cal A}u_n(t),u_n\rangle \, dt &=&\int_0^T \langle J\dot u+J{\cal A}u(t),u\rangle \, dt\\ \lim_{n\rightarrow \infty}\langle u_n(T)-u_n(0),R\frac{u_n(T)+u_n(0)}{2}\rangle
&=&
\langle u(T)-u(0),R\frac{u(T)+u(0)}{2}\rangle.
\end{eqnarray*}
3) \, Let now  ${\cal B}: D({\cal B})\subseteq L^2_X[0, T] \rightarrow L^2_X[0, T]$ be the operator defined by ${\cal B} =J\dot u(t)+J{\cal A}u(t))$ with $D({\cal B})=\big \{ u \in L^2_X: J\dot u(t)+J{\cal A}u(t)) \in L^2_X \text  {  and  } u(0)=u(T)=0\big \}$ then $R({\cal B})$ is dense in $L^2_X$.  Also, it follows from ($C'_1$) and ($C'_2$) that $\tilde H(t,.,.)$ is continuous in both variable on $L^2_X \times L^2_X$ and $\ell$ is continuous in both variables on $Y \times Y.$ Thus, for every $u\in W$, we can write
\begin{eqnarray*}
\sup\limits_{v\in W}M(u, v)&=& \sup\limits_{v\in W} \Big \{ \int_0^T\left[ \langle J\dot v (t)+J{\cal A}v(t) ,u(t)\rangle-\tilde H_L \big (t,-J\dot v-J{\cal A}v(t), -J\dot u-J{\cal A}u(t) \big ) \right]\, dt\\ & &  
+\left\langle u(T)-u(0),R\frac{v(T)+v(0)}{2}\right\rangle+\left \langle v(T)-v(0),R\frac{u(T)+u(0)}{2}\right\rangle
 \\&&-\frac{1}{\beta} \ell \big( v(T)-v(0), R\frac {v(T)+v(0)}{2}\big) \Big\}
- \int_0^T\langle J\dot u(t)+J{\cal A}u(t),u(t)\rangle \, dt
\\&&-\left\langle u(T)-u(0),R(u(T)+u(0))\right\rangle 
+\frac{1}{\beta} \ell \big( u(T)-u(0), R\frac {u(T)+u(0)}{2}\big)\\
&=& \sup\limits_{v\in W}\sup\limits_{v_0 \in D({\cal B})} \Big \{ \int_0^T\left[ \langle J\dot v (t)+J{\cal A}v(t) ,u(t)\rangle-\tilde H_L \big (t,-J\dot v-J{\cal A}v(t), -J\dot u-J{\cal A}u(t) \big ) \right]\, dt\\ & &  
+\left\langle u(T)-u(0),R\frac{(v+v_0)(T)+(v+v_0)(0)}{2}\right\rangle+\left \langle (v+v_0)(T)-(v+v_0)(0),R\frac{u(T)+u(0)}{2}\right\rangle \\&&
 -\frac{1}{\beta} \ell \big( (v+v_0)(T)-(v+v_0)(0), R\frac {(v+v_0)(T)+(v+v_0)(0)}{2}\big) \Big\}
- \int_0^T\langle J\dot u(t)+J{\cal A}u(t),u(t)\rangle \, dt\\&&-\left\langle u(T)-u(0),R(u(T)+u(0))\right\rangle 
+\frac{1}{\beta} \ell \big( u(T)-u(0), R\frac {u(T)+u(0)}{2}\big)
\end{eqnarray*}
With a change of variable $w=v+v_0$, we get
\begin{eqnarray*}
\sup\limits_{v\in W}M(u, v)
&=& \sup\limits_{w\in W} \sup\limits_{v_0 \in D({\cal B})}\Big \{ \int_0^T\left[ \langle{\cal B} w(t)-{\cal B} v_0(t) ,u(t)\rangle-\tilde H_L \big (t,-{\cal B} w(t)+{\cal B} v_0(t), -J\dot u-J{\cal A}u(t) \big ) \right]\, dt\\ & &  
+\left\langle u(T)-u(0),R\frac{w(T)+w(0)}{2}\right\rangle+\left \langle w(T)-w(0),R\frac{u(T)+u(0)}{2}\right\rangle
\\&& -\frac{1}{\beta} \ell \big( w(T)-w(0), R\frac {w(T)+w(0)}{2}\big) \Big\}
- \int_0^T\langle J\dot u(t)+J{\cal A}u(t),u(t)\rangle \, dt
\\&&-\left\langle u(T)-u(0),R(u(T)+u(0))\right\rangle 
+\frac{1}{\beta} \ell \big( u(T)-u(0), R\frac {u(T)+u(0)}{2}\big).
\end{eqnarray*}
Since $D({\cal B})$  and $R({\cal B})$ are
dense in $L^2_X$ and $x \rightarrow \int_0^T\tilde H_L(t, x(t), y(t))
 \, dt  $ is continuous in $L^2_X$ for every $y \in L^2_X$, we get that 
 \begin{eqnarray*}
\sup\limits_{v\in W}M(u, v)&=&
 \sup\limits_{w\in W} \sup\limits_{x \in L^2_X}\Big \{ \int_0^T\left[ \langle x(t) ,u(t)\rangle-\tilde H_L \big (t,-x(t), -J\dot u-J{\cal A}u(t) \big ) \right]\, dt\\ & &  
+\left\langle u(T)-u(0),R\frac{w(T)+w(0)}{2}\right\rangle+\left \langle w(T)-w(0),R\frac{u(T)+u(0)}{2}\right\rangle\\
&&-\frac{1}{\beta} \ell \big( w(T)-w(0), R\frac {w(T)+w(0)}{2}\big) \Big\}
- \int_0^T\langle J\dot u(t)+J{\cal A}u(t),u(t)\rangle \, dt
\\&&-\left\langle u(T)-u(0),R(u(T)+u(0))\right\rangle 
+\frac{1}{\beta} \ell \big( u(T)-u(0), R\frac {u(T)+u(0)}{2}\big).
\end{eqnarray*}
  Now for each $(a,b)
\in D(\tilde A)$, there is $w \in W$ such that $w(0)=a$ and $w(T)=b,$
namely the linear path $w(t)=\frac{T-t}{T}a+\frac{t}{T}b.$  Since
also $Z$ is dense in $Y$ and $\ell$ is continuous, we
finally obtain that
\begin{eqnarray*}
\sup\limits_{v\in W}M(u, v)&=&\sup\limits_{(a,b) \in Z \times Z} \sup\limits_{x \in L^2_X}\Big \{ \int_0^T\left[ \langle x(t) ,u(t)\rangle-\tilde H_L \big (t,-x(t), -J\dot u-J{\cal A}u(t) \big ) \right]\, dt\\ & &  
+\left\langle u(T)-u(0),R\frac{b+a}{2}\right\rangle+\left \langle b-a,R\frac{u(T)+u(0)}{2}\right\rangle
 -\frac{1}{\beta} \ell \big( b-a, R\frac {b+a}{2}\big) \Big\}\\&&
- \int_0^T\langle J\dot u(t)+J{\cal A}u(t),u(t)\rangle \, dt-\left\langle u(T)-u(0),R(u(T)+u(0))\right\rangle\\
&& 
+\frac{1}{\beta} \ell \big( u(T)-u(0), R\frac {u(T)+u(0)}{2}\big)\\
&=& \int_0^T L(t,u(t), -J\dot u(t)-J{\cal A}u(t)) \, dt -\int_0^T   \langle J\dot
u(t)+J{\cal A}u(t),u(t)\rangle  \, dt 
\\&&+ (\frac{1}{\beta} \ell)^* \big(R\frac{u(T)+u(0)}{2}, u(T)-u(0)\big)\\&&-\left\langle u(T)-u(0),R(u(T)+u(0))\right\rangle 
+\frac{1}{\beta} \ell \big( u(T)-u(0), R\frac {u(T)+u(0)}{2}\big)\\&=&I(u).
\end{eqnarray*}
 The following three lemmas are dedicated to the proof of the coercivity of $u\to M(u,0)$ on $W$.
\begin{lemma}\label{H-positive2}
For any $u \in W$ we have
\begin{equation*}
\ell \big( u(T)-u(0), R\frac {u(T)+u(0)}{2}\big)+
\int_0^T \langle J\dot u(t),J{\cal A}u(t)\rangle dt\geq 0.
\end{equation*}
\end{lemma}

\paragraph{Proof:} Indeed, for $u=(p,q)$ we have
\begin{eqnarray*}
\int_0^T \langle J\dot u(t),J{\cal A}u(t)\rangle dt &=&  \int_0^T  \langle (- \dot q, p),(Aq,Ap)\rangle dt\\
&=&  - \int_0^T  \langle \dot q,Aq\rangle dt+\int_0^T  \langle \dot p,Ap\rangle dt\\
&=&-\frac{1}{2} \int_0^T \frac {d}{dt}\|A^{\frac{1}{2}}q\|_Y^2\,dt+\frac{1}{2} \int_0^T \frac{d}{dt}\|A^{\frac{1}{2}}p\|_Y^2\,dt\\
&=&-\frac{1}{2}\|A^{\frac{1}{2}}q(T)\|_Y^2+\frac{1}{2}\|A^{\frac{1}{2}}q(0)\|_Y^2+\frac{1}{2}\|A^{\frac{1}{2}}p(T)\|_Y^2-
\frac{1}{2}\|A^{\frac{1}{2}}p(0)\|_Y^2\\
&=&-\langle  Aq(T)-Aq(0), \frac {q(0)+q(T)}{2}\rangle
+
\langle Ap(T)-Ap(0), \frac {p(0)+p(T)}{2}\rangle\\
&=&\langle \tilde A u(T)-\tilde A u(0), R\frac
{u(T)+u(0)}{2}\rangle\\
&\geq &- \ell \big( u(T)-u(0), R\frac {u(T)+u(0)}{2}\big).
\end{eqnarray*}
 
\begin{lemma}\label{H-estimate} For each $u \in W$, the following estimate holds:
\begin{eqnarray*}
\big | \int_0^T \langle J\dot u(t)+J{\cal A}u(t),u(t)\rangle \, dt \big | +\big |\langle u(T)-u(0),R( u(0)+u(T))\rangle \big | \leq 4c\sqrt{T}   \|u\|_W^2.
\end{eqnarray*}
\end{lemma}
\paragraph{Proof:}
It suffices to combine the following two estimates.
\begin{eqnarray*}
\big |\langle u(T)-u(0),R( u(0)+u(T))\rangle \big | &=& \big | \int_0^T\langle \dot u(t),R( u(0)+u(T))\rangle \, dt \big | \nonumber  \\ &\leq & \sqrt{T} \|\dot u\|_{L^2_X} \|u(T)+u(0)\|_{X} \leq 2\sqrt{T} \|\dot u\|_{L^2_X} \|u\|_{C(0,T;X)} \nonumber \\ & \leq & 2c\sqrt{T}   \|u\|_W^2
\end{eqnarray*}
and
\begin{eqnarray*}
\big | \int_0^T \langle J\dot u(t)+J{\cal A}u(t),u(t)\rangle \, dt \big | & \leq & \| u\|_{L^2_X} (\|\dot u\|_{L^2_X}+\|{\cal A} u\|_{L^2_X})  \nonumber \\
& \leq & \sqrt{T} \|u\|_{C(0,T;X)}  (\|\dot u\|_{L^2_X}+\|{\cal A} u\|_{L^2_X})    \nonumber\\
& \leq & 2c\sqrt{T} \|u\|_W^2.
\end{eqnarray*}
\begin{lemma}\label{H-coercive} There exists a constant $C\geq 0$ such that  for any $u \in W$:
 \begin{equation*} 
M(u, 0)\geq  (\frac{1}{2\beta}-4c\sqrt{T})\|u\|_W^2 -C. 
\end{equation*}
\end{lemma}
\paragraph{Proof:}
Note first that
\begin{eqnarray}\label{H-estimate1}
\int_0^T\tilde H_L(t,0,-J\dot u(t)-J{\cal A}u(t)) \,dt &=&\sup_{x \in L^2_X}\int_0^T\left [ \langle x(t) , -J\dot u(t)-J{\cal A}u(t)\rangle- L(t, x(t), 0)\right] \, dt \nonumber  \\
&\geq & \sup_{x \in L^2_X}\int_0^T \left  [\langle x(t) , -J\dot u(t)-J{\cal A}u(t)\rangle -\frac{\beta}{2}\| x(t)\|^2-\gamma (t)  \right ] \, dt \nonumber \\
&= & \frac {1}{2\beta} \int_0^T \|J\dot u(t)+J{\cal A}u(t)\|^2 \, dt-\int_0^T \gamma (t) \, dt \nonumber \\
&=&\frac {1}{2\beta} \int_0^T ( \|\dot u(t)\|^2+\|{\cal A}u(t)\|^2) \, dt+\frac{1}{\beta}\int_0^T \langle J\dot u(t), J{\cal A}u(t)\rangle  \, dt  \nonumber \\&&-\int_0^T \gamma (t) \, dt. 
\end{eqnarray}
It follows from Lemma \ref{H-positive2}, Lemma \ref{H-estimate} and (\ref{H-estimate1}) that  
\begin{eqnarray}\label{H-estimate2}
M(u,0) &\geq& 
\frac{1}{2\beta}\int_0^T( \|\dot u(t)\|^2+\|{\cal A}u(t)\|^2)\, dt  +\frac{1}{\beta}\int_0^T \langle J\dot u(t),J{\cal
A}u(t)\rangle\, dt-C- \int_0^T \langle J\dot u(t)+J{\cal A}u(t),u(t)\rangle\, dt \nonumber
 \\
 &&-\langle R( u(0)+u(T)), u(T)-u(0)\rangle  
+ \frac{1}{\beta} \ell \big( u(T)-u(0), R\frac {u(T)+u(0)}{2}\big) \nonumber \\
&\geq& \frac{1}{2\beta}\int_0^T( \|\dot
u(t)\|^2+\|{\cal A}u(t)\|^2)\, dt - 4c\sqrt{T}\|u\|_W^2 \nonumber\\
&& + \frac{1}{\beta}\int_0^T \langle J\dot
u(t),J{\cal A}u(t)\rangle\, dt+\frac{1}{\beta} \ell  \big(
u(T)-u(0), R\frac {u(T)+u(0)}{2}\big) -C\nonumber \\
&\geq & (\frac{1}{2\beta}-4c\sqrt{T})\|u\|_W^2+ \frac{1}{\beta} \left ( \ell  \big(
u(T)-u(0), R\frac {u(T)+u(0)}{2}\big)+\int_0^T \langle J\dot
u(t),J{\cal A}u(t)\rangle\, dt \right )-C \nonumber\\
&\geq & (\frac{1}{2\beta}-4c\sqrt{T})\|u\|_W^2-C.
\end{eqnarray}
\paragraph{Proof of Proposition \ref{Relax.2}:}
It follows from ($C'_1$) and ($C'_2$)   that ${\cal L}$ is finite on $W \times W$, and from Lemma \ref{H-positive1}  that  $I$ is selfdual on $W$. In view of the coercivity guaranteed by Lemma \ref{H-coercive}, we can apply Theorem \ref{two} to get $v\in W$ such that $I(v)=0$. It follows that 
\begin{eqnarray*}
L(t,v(t),-J\dot v(t)-J{\cal A}v(t))-\langle J\dot
v(t)+J{\cal A}v(t),v(t)\rangle = 0 \quad  \quad \text{for }t \in [0,T].
\end{eqnarray*}
and 
\begin{eqnarray*}
\frac {1}{\beta}\ell \big( v(T)-v(0), R\frac {v(T)+v(0)}{2}\big) +(\frac {1}{\beta}\ell)^*\big( R\frac {v(T)+v(0)}{2}, v(T)-v(0)\big) -\langle v(T)-v(0),R(v(0)+v(T))\rangle = 0.
\end{eqnarray*}
and we are done with the proposition.
\hfill $\square$

\paragraph{Proof of Theorem \ref{thm.new.var}:}
We  just need to show that  the result of Proposition \ref{Relax.2} still holds if one  replaces  ($C'_1$) and ($C'_2$)  with ($C_1$)  and ($C_2$) respectively. Indeed, for $0<\la < \frac {1}{8 c\sqrt{T}}- \beta$,   we replace $L$ with $L^2_{\la}$ in such a way that 
\begin{eqnarray}\label{la-estimate1}
L^2_{\la}(x,p)=\inf \{L(x,r)+ \frac {\|p-r\|^2}{2 \la}+ \frac{\la \|x\|^2}{2}; r \in X \} & \leq & L(x,0)+ \frac {\|p\|^2}{2 \la}+ \frac{\la \|x\|^2}{2} \nonumber \\
&\leq &   \frac{\la+\beta}{2} \|x\|^2+\frac {\|p\|^2}{2 \la}+\gamma(t)
\end{eqnarray}
and therefore satisfies ($C'_1$) whenever $\la+\beta <\frac {1}{8 c\sqrt{T}}.$ 
We also replace $\ell$ with the Lagrangian defined on $Y\times Y$ as
\begin{eqnarray*}
\ell^{1,2}_{\lambda}(x,r)=\inf \big \{   \ell(y,s)+ \frac{1}{2 \la}\|x-y\|_Y^2+  \frac{\la}{2 }\|r\|_Y^2+ \frac{1}{2 \la}\|s-r\|_Y^2+  \frac{\la}{2 }\|y\|^2; \, y \in Y, s \in Y\big \}
\end{eqnarray*}
and by $+\infty$ if either $x$ or $p$ belonds to $X\setminus Y$.  It is easily seen that $\ell^{1,2}_{\la}$ satisfies ($C'_2$) on $Y\times Y$. Moreover, $\ell^{1,2}_{\la}(x,p)\geq -\langle x, \tilde Ap \rangle $ for all $p \in Y$. Indeed, it clearly suffices to assume that $x\in Y$ also. Since $\ell$ is lower semi-continuous on $Y\times Y$,  there is $x_{\la}, p_{\la} \in Y$ such that  
\begin{eqnarray*}
{\ell}^{1,2}_{\lambda}(x,p)= \ell(x_{\la},p_{\la})+ \frac{1}{2 \la}\|x-x_{\la}\|_Y^2+  \frac{\la}{2 }\|p\|_Y^2+ \frac{1}{2 \la}\|p_{\la}-p\|_Y^2+  \frac{\la}{2 }\|x_{\la}\|_Y^2.
\end{eqnarray*}
It follows that
\begin{eqnarray*}
{\ell}^{1,2}_{\lambda}(x,p)&\geq & -\langle x_{\la}, \tilde A p_{\la} \rangle -  \langle x-x_{\la}, p \rangle_{Y \times Y} -\langle x_{\la},p- p_{\la} \rangle_{Y \times Y}\\ 
&=&  -\langle x_{\la}, p_{\la} \rangle_{Y \times Y} -  \langle x-x_{\la}, p \rangle_{Y \times Y} -\langle x_{\la},p- p_{\la} \rangle_{Y \times Y}\\&=& -\langle x, p \rangle_{Y \times Y}=-\langle x, \tilde Ap \rangle.
\end{eqnarray*}
We can now apply Proposition \ref{Relax.2}, to find $u_{\la} \in W$ with 
\begin{eqnarray}\label{la-estimate2}
I_{\la}(u_{\la})&=& \int_0^T \left\{ L^2_{\la} (t,u_{\la}(t),-J\dot u_{\la} (t)-J{\cal A}u_{\la}(t) )-\langle J \dot u_{\la}(t)+J{\cal A}u_{\la}(t),u_{\la}(t)\rangle \right\} dt \nonumber\\
 &&+ \frac{1}{\beta+\la} \ell^{1,2}_{\la} \big( u_{\la}(T)-u_{\la}(0), R\frac {u_{\la}(T)+u_{\la}(0)}{2}\big) \nonumber\\&&+( \frac{1}{\beta+\la} \ell^{1,2}_{\la} )^*\big( R\frac {u_{\la}(T)+u_{\la}(0)}{2}, u_{\la}(T)-u_{\la}(0)\big)\nonumber\\
&& -\langle u_{\la}(T)-u_{\la}(0),R( u_{\la}(0)+u_{\la}(T))\rangle=0.
\end{eqnarray}
It follows from  (\ref{la-estimate1}) and part (1) of Lemma \ref{mild.estimate} that
\begin{eqnarray}\label{la-estimate2}
\int_0^T L^2_{\la} (t,u_{\la}(t),-J\dot u_{\la} (t)-J{\cal A}u_{\la}(t) )\, dt \geq \frac{1}{2(\la + \beta)}\|J\dot u_{\la} (t)+J{\cal A}u_{\la}(t) \|^2_{L^2_X}-C_2.
\end{eqnarray}
From (\ref{la-estimate2}),  (\ref{la-estimate1}),    Lemma \ref{H-estimate} and the fact that $( \frac{1}{\beta+\la} \ell^{1,2}_{\la} )^*$ is bounded from below, we get  that 
\begin{eqnarray}
\frac{1}{2(\la + \beta)}\|J\dot u_{\la} (t)+J{\cal A}u_{\la}(t) \|^2_{L^2_X}-4c\sqrt{T}\|u_{\la}\|_W^2 + \frac{1}{\beta+\la} \ell^{1,2}_{\la} \big( u_{\la}(T)-u_{\la}(0), R\frac {u_{\la}(T)+u_{\la}(0)}{2}\big)\leq C
\end{eqnarray}
where $C$ is a constant independent of $\la.$ By the same argument as in (\ref{H-estimate2}) we obtain
\begin{eqnarray}
(\frac{1}{2(\la + \beta)}-4c\sqrt{T})\|u_{\la}\|_W^2\leq C, 
\end{eqnarray}
which ensures the boundedness of $u_{\la}$ in $W$.  Assuming $u_{\la} \rightharpoonup u$ weakly in $W$, it follows from Lemmas \ref{2-reg} 
that
$I(u)\leq \liminf_{\la}I_{\la}(u_{\la})=0.$
Since on the other hand $I(u)\geq 0$, the latter is therefore equal zero and  $u$ is a solution of (\ref{H-equation01}). \hfill $\square$

\subsection{Coercive Hamiltonian systems of PDEs}

We shall now apply Theorem \ref{thm.new.var} to the ASD-Lagrangian  
$L(t, u, p)=\phi (t, u) + \phi^*(t, -J{\cal B}u -p)$  on $X\times X$, where $\phi:[0,T]\times X \ra \R $ is a time-dependent  convex lower semi-continuous function on $X$.

As to the boundary Lagrangian, we shall associate to a given convex lower semi-continuous function $\psi$ on $X$, the following function on $X$
\begin{equation}
\psi^o(p)=\sup\{\langle p, \tilde Ax\rangle -\psi (x);\, x\in Y\}.
\end{equation} 
It is clear that $\ell$ is convex and lower semi-continuous on $X$, and that the function $\ell (x,p)=\psi (x) +\psi^o(-p)$ satisfies
\begin{equation}
\hbox{$\ell (x,p) \geq -\langle \tilde Ax, p\rangle$ for all $x\in Y, p\in X$.} 
\end{equation}
 It is also easy to see that 
 \begin{equation}
 \hbox{$\ell^*(p,x)=\psi^*(p)+\psi(-\tilde A^{-1}x)$ for all $x, p\in X$.} 
  \end{equation}
We can obtain the following. 
\begin{theorem} \label{Ham-phi-psi.1} Let $A$ be a linear operator on $H$ as above and let $\cal B$ be an operator on $X$ such that $J{\cal B}$ is skew adjoint. Let  $ \psi$ be a convex lower semi-continuous function on $X$ that is bounded below and such that $0\in {\rm Dom}(\psi)$, and let $\phi:[0,T]\times X \ra \R $ be a time-dependent  convex lower semi-continuous function on $X$ satisfying for some $\beta>0$,  $\gamma , \al  \in L^2 (0,T;\mathbb{R_{+}})$
\begin{equation}\label{b1}
\hbox{$- \al (t)\leq \phi(t,u)+\phi^*(t, J{\cal B} u)\leq\frac{\beta}{2} \|u\|_X^2+\gamma (t)$ for every $u\in X$ and a.e. $t\in [0,T]$.}
\end{equation}
Assume that 
\begin{equation} 0<T< \frac{1}{64c^2\beta^2},
\end{equation}
then  the infimum on $W$ of the functional
\begin{eqnarray} \label{new.var2.phi}
  I(u)&=& \int_0^T \left\{ \phi \big (t, u(t) \big)+ \phi^*\big (t, J\dot u (t)+J{\cal A}u(t)+J{\cal B}u(t) \big )-\langle J\dot u(t)+J{\cal A}u(t),u(t)\rangle \right\} dt \nonumber\\
&&+ \frac{1}{\beta} \psi \big( u(T)-u(0)) + \frac{1}{\beta} \psi^o \big( R\frac {u(T)+u(0)}{2}\big)\nonumber\\   
 &&+\frac{1}{\beta} \psi^* \big( \beta R\frac {u(T)+u(0)}{2}\big)+\frac{1}{\beta} \psi \big( \beta \tilde A^{-1}( u(T)-u(0))\big)\nonumber \\   
&& -\langle u(T)-u(0),R( u(0)+u(T))\rangle
\end{eqnarray}
 is  equal to zero and is attained at some $v \in W $ which is then a solution of the following system:
\begin{eqnarray*}\left\{ \begin{array}{lcl}
J\dot v(t)+J{\cal A}v(t)+J{\cal B}v(t) &= & \partial \phi \big( t,v(t)\big) \quad \hbox{\rm a.e on $[0, T]$}\\
\hfill R\frac {v(T)+v(0)}{2} &\in & \frac {1}{\beta}\partial
\psi \big(v(T)-v(0) \big).
\end{array}\right.
\end{eqnarray*}
\end{theorem}
{\bf Proof:} We apply Theorem \ref{thm.new.var} to the ASD-Lagrangian  
$L(t, u, p)=\phi (t, u) + \phi^*(t, -J{\cal B}u -p)$  on $X\times X$, where $\phi:[0,T]\times X \ra \R $ is a time-dependent  convex lower semi-continuous function on $X$.\\
As to the boundary Lagrangian, consider a convex lower semi-continuous function $\psi$ on $X$ that is bounded below and such that $0\in {\rm Dom}(\psi)$, and define on $X\times X$, the convex function 
$\ell (x,p)=\psi (x) +\psi^o(-p)$.\\
Suppose now that for some $x, y\in Y$, we have $\ell(x, y)+\ell^*(y, x)-2\langle x, y\rangle=0$. 
This means that 
\[
\psi (x)+\psi^o(-y) +\psi^*(y)+\psi(-\tilde A^{-1}x)-2\langle x, y\rangle=0. 
\]
Since $\psi^o(-y) +\psi(-\tilde A^{-1}x) \geq \langle x, y\rangle$, it follows that $\psi (x)+\psi^*(y)=\langle x, y\rangle$ from which we conclude that $y\in \partial \psi (x)$.  
 
\begin{remark} Here again, the general boundary conditions we obtain will allow us to obtain periodic and other type of solutions.  Indeed, 
  \begin{itemize}
\item Periodic solutions $v(0)=v(T)$, then $\psi$ is chosen as:
\begin{eqnarray*}\psi(w)=\left\{\begin{array}{ll}
0 \quad &w=0\\
+\infty &\mbox{elsewhere}.\end{array}\right.
\end{eqnarray*}
\item Anti periodic solutions $v(0)=-v(T)$, then $\psi \equiv0.$
\item Initial boundary condition $p(0)=p_0$ and $q(T)=q_0$ for a given $p_0, q_0\in
H$. Let $v_0=(-p_0,q_0)$ and $\psi (w)= \frac{\beta}{4}\|w\|^2-\beta \braket{w}{v_0}$. then it follows 
that
$$R\frac {v(T)+v(0)}{2}=\frac{1}{\beta} \partial \big [\psi (v(T)-v(0))\big ]
= \frac{v(T)-v(0)}{2}-{v_0}.$$
Setting $v=(p,q)$ we have
\begin{eqnarray*}
 (\frac {p(T)+p(0)}{2},-\frac {q(T)+q(0)}{2})&=& 
 R (\frac {p(T)+p(0)}{2},\frac {q(T)+q(0)}{2})\\
&=&R\frac {v(T)+v(0)}{2}  \\
&=& \frac{v(T)-v(0)}{2}-{v_0} \\
&=&  (\frac{p(T)-p(0)}{2}+p_0,\frac{q(T)-q(0)}{2}-q_0)
\end{eqnarray*}
from which we obtain $p(0)=p_0$ and $q(T)=q_0.$
\end{itemize}
 \end{remark}

\subsubsection*{Example 6:  A  coercive Hamiltonian System  involving the bi-Laplacian}
Let $\Omega$ be a bounded domain in $\R^N$ and consider the following Hamiltonian System,
\begin{eqnarray}\label{ham-ex20}\left\{ \begin{array}{lcl}
-\dot v(t)+\Delta^2 v - \Delta v&= & \partial \phi_1 (t,u) \qquad (t,x) \in (0,T) \times \Omega,\\
\hfill  \dot u(t)+\Delta^2 u+ \Delta u&= &\partial \phi_2(t,v)  \qquad (t,x) \in (0,T) \times \Omega,\\
\hfill{}u=\Delta u &=& 0 \hfill{}    (t,x) \in [0,T] \times  \partial \Omega\\
\hfill{}v=\Delta v &=& 0 \hfill{}  (t,x) \in [0,T] \times  \partial \Omega
\end{array}\right.
\end{eqnarray}
where $\phi_i, i=1,2$ are two convex lower semi-continuous functions on $H:=H_0^1(\Omega)$ considered as a  Hilbert space with the inner product $\langle u, v \rangle =\int_{\Omega}\nabla u \cdot \nabla v\, dx$. We consider $Y=\{u \in H_0^1(\Omega); \Delta u \in H_0^1(\Omega)  \}$ equipped with the norm
$\|u\|^2_Y= \int_{\Omega} |\nabla \Delta u|^2 \, dx$. Theorem \ref{Ham-phi-psi.1} yields the following existence result. 
\begin{theorem} Suppose  $\phi_1$ and $\phi_2$ satisfy the following condition:
\begin{eqnarray}\label{ham-ex20-1}
\gamma_i (t)+ c_i \|u\|^2_{L^2(\Omega)} \leq \phi_i (t,u) \leq   \alpha_i (t)+ C_i \|u\|^2_{H^1_0(\Omega)}  \qquad i=1,2
\end{eqnarray}
where  $ \gamma_i ,  \alpha_i  \in L^2([0, T]),$ and $c_i,C_i>0.$  Then for $T$ small enough  there exist $u,v \in W$ satisfying (\ref{ham-ex20})
with either of the following boundary conditions,
 \begin{itemize}
\item Periodic solutions $u(0)=u(T)$ and $v(0)=v(T).$
\item Anti periodic solutions $u(0)=-u(T)$ and $v(0)=-v(T).$
\item Initial boundary condition $u(0)=u_0$ and $v(T)=v_0$ for a given $v_0, u_0\in H.$
\end{itemize}
\end{theorem}
\paragraph{Proof} Let $Au=\Delta^2u$ so that  for $U=(u,v)$, ${\cal A}U={\cal A}(u,v)=(\Delta^2u,-\Delta^2 v)$. Consider the operator ${\cal B}U=(\Delta u,  \Delta v)$ in such a way that $J{\cal B}U=(-\Delta v, \Delta u)$ is skew-adjoint on $H^1_0(\Omega)\times H^1_0(\Omega)$. Equation (\ref{ham-ex20}) can be rewritten as follows 
\begin{eqnarray*}
J \dot U (t)+J{\cal A}U(t)= \bar \partial L(t, U(t))
\end{eqnarray*}
where $L(t,U,V)= \Phi (t,U)+\Phi^*(t,JBU-V)$ with $\Phi(t,U)=\phi_1 (t,u)+\phi_2(t,v).$ We just need to show that $L$ satisfies condition $(C_1)$ in Theorem \ref{thm.new.var}.  Let $C=\max\{C_1,C_2\}$, $c=\min\{c_1,c_2\}$, $\gamma (t)= \min \{\gamma_1(t), \gamma_2(t)\}$ and $ \alpha(t)=\max\{\alpha_1(t), \alpha_2(t)\}.$ It follows from (\ref{ham-ex20-1}) that
\begin{eqnarray*}
\gamma (t)+ c\|U\|^2_{L^2(\Omega)} \leq \Phi (t,U) \leq   \alpha (t)+ C \|U\|^2_{H^1_0(\Omega)}, 
\end{eqnarray*}
and therefore 
\begin{eqnarray*}
 -\alpha (t)+ \frac{1}{4C} \|U\|^2_{H^1_0(\Omega)} \leq \Phi^*(t,U) \leq   -\gamma (t)+ \frac{1}{4c}\|\nabla (-\Delta)^{-1} U\|^2_{L^2(\Omega)},
\end{eqnarray*}
from which we obtain 
\begin{eqnarray*}
\gamma (t) -\alpha (t)\leq   L(t,U,0)  &\leq & \alpha (t)-\gamma (t)+ C \|U\|^2_{H^1_0(\Omega)}  + \frac{1}{4c}\|\nabla (-\Delta)^{-1} J{\cal B}U\|^2_{L^2(\Omega)}\\
&=&\alpha (t)-\gamma (t)+C \|U\|^2_{H^1_0(\Omega)}  + \frac{1}{4c}\|\nabla U\|^2_{L^2(\Omega)}\\
&=&\alpha (t)-\gamma (t)+(C+\frac{1}{4c}) \|U\|^2_{H^1_0(\Omega)}.
\end{eqnarray*}
Hence for $T$ small enough,  Theorem \ref{Ham-phi-psi.1} applies to yield our claim. 

\subsection{Non-coercive Hamiltonian systems of PDEs}
Under a certain commutation property, we can relax the boundedness condition (\ref{b1}) provided one settles for periodic solutions up to an isometry. 

\begin{corollary} \label{cor.new.var}
Let $L:[0,T]\times X \times X \ra \R \cup \{+\infty\}$, $\ell: X \times X \ra \R \cup \{+\infty\}$,  and $A: D(A)\subset H\to H$ be as in Theorem \ref{thm.new.var}, and let ${\cal B}$ be a skew-adjoint operator on $H\times H$ such that ${\cal A}{\cal B}={\cal B}{\cal A}$ on $D({\cal A})$, and let $(S_t)_t$ be its corresponding  $C_0$-unitary group of operators on $X$. 
Then  the infimum of the functional
\begin{eqnarray} \label{new.var2}
 I(u)&=& \int_0^T \left\{ L (t,S_tu(t),-JS_t\dot u (t)-J{\cal A}S_tu(t) )-\langle JS_t \dot u(t)+J{\cal A}S_tu(t),S_tu(t)\rangle \right\} dt \nonumber\\
 &&+ \frac{1}{\beta} \ell \big( u(T)-u(0), R\frac {u(T)+u(0)}{2}\big) +( \frac{1}{\beta} \ell )^*\big( R\frac {u(T)+u(0)}{2}, u(T)-u(0)\big)\nonumber\\
&& -\langle u(T)-u(0),R( u(0)+u(T))\rangle
\end{eqnarray}
on $W$
 is  equal to zero and is attained at some $u \in W $ in such a way that  $v(t):=S_tu(t)$ is  a solution of
\begin{eqnarray}\label{H-equation2}\left\{ \begin{array}{lcl}
\hfill  J\dot v(t)+J{\cal A}v(t)+J{\cal B}v(t) &= &\bar \partial L\big( t,v(t)\big)\\
\hfill R\frac {S_{(-T)}v(T)+v(0)}{2}&= & \bar \partial \frac {\ell}{\beta} \big(S_{(-T)}v(T)-v(0) \big).
\end{array}\right.
\end{eqnarray}
\end{corollary}
\paragraph{Proof:} 
It follows from Proposition \ref{permanence} that $L_S(t,x,y):=L(t, S_tx, S_ty)$ is anti-self dual Lagrangian on $[0,T] \times X \times X.$ Since $S_t$ is norm preserving, assumption ($C_1$ ) holds for the new Lagrangian $L_S$.   Therefore there exists $u \in W$ such that $I(u)=0$  and  $u$ is a solution of
\begin{eqnarray} \label{equ-semigroup0}\left\{ \begin{array}{lcl}
\hfill  J\dot u(t)+J{\cal A}u(t) &= &\bar \partial L_S \big( t,u(t)\big)\\
\hfill R\frac {S_{(-T)}u(T)+u(0)}{2}&= & \bar \partial \frac {\ell}{\beta} \big(S_{(-T)}u(T)-u(0) \big).
\end{array}\right.
\end{eqnarray}
Note that 
$\bar \partial L_S \big( t,u(t)\big)=S^*_t\bar \partial L (t, S_t u(t))$
which  together with equation (\ref{equ-semigroup0}),  imply that 
\begin{eqnarray*}
S_t \big (  J\dot u(t)+J{\cal A}u(t) \big )=\bar \partial L (t, S_t u(t)).
\end{eqnarray*}
Since ${\cal A}{\cal B}={\cal B}{\cal A}$ on $D({\cal A})$, we have $S_t{\cal A}u(t)={\cal A}S_tu(t)$ and therefore 
\begin{eqnarray} \label{equ-semigroup}\left\{ \begin{array}{lcl}
\hfill  JS_t \dot u(t)+J{\cal A}S_t u(t) &=&\bar \partial L (t, S_t u(t))\\
\hfill R\frac {u(T)+u(0)}{2}&= & \bar \partial \frac{\ell}{\beta}\big(u(T)-u(0), R\frac {u(T)+u(0)}{2} \big).
\end{array}\right.
\end{eqnarray}
To show that $v(t):=S(t)u(t)$ is a solution of problem (\ref{H-equation2}), 
substitute $u(t)=S(-t)v(t)$ in (\ref{equ-semigroup}) to obtain
\begin{eqnarray}\left\{ \begin{array}{lcl}
\hfill  J\dot v(t)+J{\cal A}v(t)+J{\cal B}v(t) &= &\bar \partial L\big( t,v(t)\big)\\
\hfill R\frac {S_{(-T)}v(T)+v(0)}{2}&= & \bar \partial \frac {\ell}{\beta} \big(S_{(-T)}v(T)-v(0)\big).
\end{array}\right.
\end{eqnarray}
\hfill $\square$\\
By applying again the above to the ASD-Lagrangian  
$L(t, u, p)=\phi (t, u) + \phi^*(t, -p)$  on $X\times X$, and $\ell (x,p)=\psi (x)+\psi^o(-p)$, 
we get the following.
\begin{theorem} \label{Ham-phi-psi.2} Let $L:[0,T]\times X \times X \ra \R \cup \{+\infty\}$, $\ell: X \times X \ra \R \cup \{+\infty\}$,  and $A: D(A)\subset H\to H$ be as in Theorem \ref{thm.new.var},  and let ${\cal B}$ be a skew-adjoint operator on $H\times H$ such that ${\cal A}{\cal B}={\cal B}{\cal A}$ on $D({\cal A})$, and let $(S_t)_t$ be its corresponding  $C_0$-unitary group of operators on $X$. Let  $ \psi$ be a convex lower semi-continuous function on $X$ that is bounded below and such that $0\in {\rm Dom}(\psi)$, and let $\phi:[0,T]\times X \ra \R $ be a time-dependent G\^ateaux differentiable  convex function on $X$ satisfying for some $\beta>0$,  $\gamma , \al  \in L^2 (0,T;\mathbb{R_{+}})$
\begin{equation} \label{b2}
\hbox{$- \al (t)\leq \phi(t,u)\leq\frac{\beta}{2} \|u\|_X^2+\gamma (t)$ for every $u\in X$ and a.e. $t\in [0,T]$.}
\end{equation}
Assuming that 
\begin{equation} 0<T< \frac{1}{64c^2\beta^2},
\end{equation}
then  the infimum on $W$ of the functional
\begin{eqnarray} \label{new.var2.phi}
\bar  I(u)&=& \int_0^T \left\{ \phi \big (t,S_tu(t) \big)+ \phi^*\big (JS_t\dot u (t)+J{\cal A}S_tu(t) \big )-\langle JS_t \dot u(t)+J{\cal A}S_tu(t),S_tu(t)\rangle \right\} dt \nonumber\\
 &&+ \frac{1}{\beta} \ell \big( u(T)-u(0), R\frac {u(T)+u(0)}{2}\big) 
+( \frac{1}{\beta} \ell)^*\big( R\frac {u(T)+u(0)}{2}, u(0)-u(T)\big)\nonumber\\
&& -\langle u(T)-u(0),R( u(0)+u(T))\rangle
\end{eqnarray}
 is  equal to zero and is attained at some $u \in W $ in such a way that  $v(t):=S_tu(t)$ is  a solution of the foillowing system:
\begin{eqnarray*}\left\{ \begin{array}{lcl}
J\dot v(t)+J{\cal A}v(t)+JBv(t) &= & \partial \phi \big( t,v(t)\big) \quad \hbox{\rm a.e on $[0, T]$}\\
\hfill R\frac {S_{(-T)}v(T)+v(0)}{2} &\in & \frac {1}{\beta}\partial
\psi \big(S_{(-T)}v(T)-v(0) \big).
\end{array}\right.
\end{eqnarray*}
\end{theorem}
 
  \subsubsection*{Example 7: Periodic solutions up to an isometry for a noncoercive Hamiltonian system involving the bi-Laplacian}
We now consider the following Hamiltonian System,
\begin{eqnarray}\label{ham-ex2}\left\{ \begin{array}{lcl}
-\dot v(t)+\Delta^2 v - \Delta u&= & \partial \phi_1 (t,u) \qquad (t,x) \in (0,T) \times \Omega,\\
\hfill  \dot u(t)+\Delta^2 u- \Delta v&= &\partial \phi_2(t,v)  \qquad (t,x) \in (0,T) \times \Omega,\\
\hfill{}u=\Delta u &=& 0 \hfill{}    (t,x) \in [0,T] \times  \partial \Omega\\
\hfill{}v=\Delta v &=& 0 \hfill{}  (t,x) \in [0,T] \times  \partial \Omega
\end{array}\right.
\end{eqnarray}
where again $\phi_i, i=1,2$ are two convex lower semi-continuous functions on $H:=H_0^1(\Omega)$ considered as a  Hilbert space with the inner product $\langle u, v \rangle =\int_{\Omega}\nabla u \cdot \nabla v\, dx$. We consider $Y=\{u \in H_0^1(\Omega); \Delta u \in H_0^1(\Omega)  \}$ equipped with the norm
$\|u\|^2_Y= \int_{\Omega} |\nabla \Delta u|^2 \, dx$. Theorem \ref{Ham-phi-psi.2} yields the following existence result. 
\begin{theorem} Suppose  $\phi_1$ and $\phi_2$ satisfy the following condition:
\begin{eqnarray}\label{ham-ex2-1}
\gamma_i (t) \leq \phi_i (t,u) \leq   \alpha_i (t)+ C_i \|u\|^2_{H^1_0(\Omega)}  \qquad i=1,2
\end{eqnarray}
where  $ \gamma_i ,  \alpha_i  \in L^2([0, T]),$ and $c_i,C_i>0.$  Then for $T$ small enough,   there exist $u,v \in W$ satisfying (\ref{ham-ex2})
with either of the following boundary conditions
  \begin{itemize}
\item Periodic solutions up to an isometry.
\item Anti periodic solutions up to an isometry.
\item Initial boundary condition $u(0)=u_0$ and $v(T)=v_0$ for a given $v_0, u_0\in H.$
\end{itemize}
\end{theorem}
\paragraph{Proof} Let again $Au=\Delta^2u$ in such a way that  for $U=(u,v)$, ${\cal A}U={\cal A}(u,v)=(\Delta^2u,-\Delta^2 v)$. Consider however the skew adjoint operator ${\cal B}U=(-\Delta v,  \Delta u)$ in such a way that $J{\cal B}U=(-\Delta u, -\Delta v)$. Problem (\ref{ham-ex2}) can be rewritten as 
\begin{eqnarray} 
 J\dot v(t)+J{\cal A}v(t)+J{\cal B}v(t) = \bar \partial L \big( t,v(t)\big)
\end{eqnarray} 
where $L(t,U,V)= \Phi (t,U)+\Phi^*(t,-V)$ with $\Phi(t,U)=\phi_1 (t,u)+\phi_2(t,v)$. In order  to show that $L$ satisfies condition $(C_1)$ in Theorem \ref{thm.new.var}, it suffices to notice that 
\begin{eqnarray*}
\gamma (t) -\alpha (t)&\leq&   L(t,U,0)=\Phi (t,U)+\Phi^*(t,0)\\ &\leq & \alpha (t)-\gamma (t)+ C \|U\|^2_{H^1_0(\Omega)} 
\end{eqnarray*}
where again  $C=\max\{C_1,C_2\}$, $\gamma (t)= \min \{\gamma_1(t), \gamma_2(t)\}$ and $ \alpha(t)=\max\{\alpha_1(t), \alpha_2(t)\}$.

\subsection*{Example 8: Periodic solutions up to an isometry for a noncoercive Hamiltonian System involving the Laplacian and transport} 
Consider the following Hamiltonian system of PDEs:
\begin{eqnarray}\label{ham-ex1}\left\{ \begin{array}{lcl}
-\dot v(t)-\Delta (v+u)+b.\nabla v&= &|u|^{p-2}u+ g(t,x)  \qquad (t,x) \in (0,T) \times \Omega,\\
\hfill \dot u(t)-\Delta (u+v)+ a.\nabla u&= &|v|^{q-2}v+ f(t,x) \qquad (t,x) \in (0,T) \times \Omega,
\end{array}\right.
\end{eqnarray}
where $a,b \in \R^N$ are two constant vectors. Let $H=L^2(\Omega)$ and $Y=H_0^1(\Omega).$  
\begin{theorem} Suppose $f,g \in L^2_H$ and $1<p,q <2.$ Then for any $T>0$ there exists $u,v \in W$ satisfying (\ref{ham-ex1}) with either of the following boundary conditions
  \begin{itemize}
\item Periodic solutions up to an isometry.
\item Anti periodic solutions up to an isometry.
\item Initial boundary condition $u(0)=u_0$ and $v(T)=v_0$ for a given $v_0, u_0\in H.$
\end{itemize}
\end{theorem}
\paragraph{Proof} Problem (\ref{ham-ex1}) can be rewritten as
\begin{eqnarray} 
J\dot U(t)+J{\cal A} U(t)+JBU(t)= \bar \partial L(t,U(t))
 \end{eqnarray}
where $ {\cal A}(u,v)=(- \Delta u,  \Delta v)$, ${\cal B}(u,v)=(-\Delta v+a.\nabla u, \Delta u-b.\nabla v)$ and $L(t,U,V)=\Phi (t,U)+\Phi^* (t,-V)$ with
\[\Phi(t,U)= \frac{1}{p}\int_{\Omega}|u|^p \,dx+\langle u, f(t,x) \rangle + \frac{1}{q}\int_{\Omega}|v|^q \,dx+\langle v, g(t,x) \rangle\]  
It is clear that all hypothesis of Theorem \ref{cor.new.var}  are satisfied.

  \section{Schr\"odinger and other nonlinear evolutions}
  
  Considering again that $X \subseteq H \subseteq X^*$ is  an evolution triple, we shall denote by $D$ the duality map between $X$ and $X^*$. We need the following  notion which is the analogue of the Palais-Smale condition (\cite{Ek.book} \cite{St}) for selfdual variational calculus.

 \begin{definition} \rm Let $L$ be  a time-dependent anti-selfdual Lagrangian on $[0,T] \times X \times X^*$, $\ell$ an anti-selfdual Lagrangian on $H\times H$, and let $\Lambda:{\cal X}_{p,q}\to L^q_{X^{*}}$ be a given map. 
 Say that $(L, \ell) $ is {\it $\Lambda$-coercive} 
 if  any  sequence $\{x_n\}_{n=1}^{\infty} \subseteq {\cal X}_{p,q}$ satisfying 
 \begin{eqnarray}\label{principle3}
 \left \{ \begin{array}{lcl}
  \hfill \dot {x}_n(t) +\Lambda x_n(t)-\frac{1}{n}\|u_n\|^{p-2}D u_n& = & -\bar \partial L (t,x_n(t)), \\ \label{eqn:zero2}
  \hfill  \frac{v_n(0)+v_n(T)}{2}
 & \in & - \bar \partial \ell\big(v_n(0)-v_n(T))
\end{array}\right.
  \end{eqnarray}
  is bounded in ${\cal X}_{p,q}$.
 \end{definition}
   The following variational principle for nonlinear evolutions established in \cite{GM3} already allows us to deal with certain Schr\"odinger equations.
  
\begin{theorem}  \label{nonlinear} {\rm  \cite{GM3}} Let $X\subset H\subset X^*$  be an evolution triple where $X$ is a reflexive Banach space, and $H$ is a Hilbert space. For $p>1$ and $q=\frac{p}{p-1}$, assume that $\Lambda:{\cal X}_{p,q}\to  L^q_{X^{*}}$  is  a  regular map such that for some nondecreasing continuous real function $w$, and $0\leq k<1$, it satisfies
 \begin{equation} \label{Lambda-10}
\hbox{$\|\Lambda x\|_{L^q_{X^*}}\leq k\|\dot x\|_{L^q_{X^*}}+ w(\|x\|_{L^p_{X}})$ for every $x\in {\cal X}_{p,q}$,}
\end{equation}
and
\begin{eqnarray}\label{Lambda-20}
\hbox{$\big |\int_0^T \langle \Lambda  x(t),  x(t)\rangle \, dt  \big | \leq  w(\|x\|_{L^p_X})$ for every $x\in {\cal X}_{p,q}$}.
\end{eqnarray}
 Let $\ell$ be an anti-selfdual Lagrangian on $H \times H$ that is bounded below with  $0 \in {\rm Dom}(\ell)$, and  let $L$ be a time dependent  anti-selfdual Lagrangian on  $[0,T] \times X \times X^*$  such  that for some $C>0$ and $r>1$, we have 
\begin{equation}\label{L.C}
\hbox{$-C \leq \int_0^TL(t, u(t), 0) dt \leq C(1+\|u\|^r_{L^p_X})$ for every $u\in L^p_X$.}
\end{equation} 
If $(L, \ell)$ is $\Lambda$-coercive, then the following functional 
\begin{equation}
I(u)= \int_0^T \Big [  L  (t,  u(t),\dot { u}(t)+\Lambda  u(t)) +\langle \Lambda  u(t),  u(t) \rangle \Big ] \, dt +\ell (u(0)- u(T), \frac { u(T)+ u(0)}{2})
\end{equation}
attains its minimum at $v \in {\cal X}_{p,q}$ in such a way that $I(v)=\inf_{u\in  {\cal X}_{p,q}}I(u)=0$ and
\begin{eqnarray}\label{principle3}
 \left \{ \begin{array}{lcl}
  \hfill \dot {v}(t)+\Lambda v(t)& = & -\bar \partial L (t,v(t)), \\ \label{eqn:zero2}
  \hfill  \frac{v(0)+v(T)}{2}
 & \in & -\bar \partial \ell\big(v(0)-v(T)).
\end{array}\right.
  \end{eqnarray}
\end{theorem}

\subsection{Initial-value Schr\"odinger evolutions}
Consider the following nonlinear Schr\"odinger equation
\begin{eqnarray}\label{sch-equ.1}
iu_t+\Delta u-|u|^{r-1}u= -i\bar \partial L(t,u) \qquad (t,x)\in [0,T] \times \Omega,
\end{eqnarray}
where $\Omega$ is a bounded domain in $\mathbb{R}^N$, and $L$ is a time dependent anti-selfdual Lagrangian on $[0, T]\times H_0^1(\Omega) \times H^{-1}(\Omega)$.  Equation (\ref{sch-equ.1}) can be rewritten as 
\begin{eqnarray*} 
u_t+\Lambda u= -\bar \partial L (t,u) \qquad (t,x)\in [0,T] \times \Omega,
\end{eqnarray*}
where $\Lambda u=-i \Delta+i|u|^{r-1}u$.  We can then deduce the following existence.
\begin{theorem} Suppose $1\leq r\leq \frac{N}{N-2}$. Let  $p=2r$  and  assume that  $L$ satisfies  
\begin{equation}\label{L.C.1}
\hbox{$-C \leq \int_0^TL(t, u(t), 0) dt \leq C(1+\|u\|^r_{L^p_{H^1_0}})$ for every $u\in L^p_{H^1_0}[0, T]$.}
\end{equation} 
\begin{eqnarray}\label{pos}
\hbox{$\langle \bar \partial L (u),-\Delta u+|u|^{r-1}u \rangle \geq 0$ for each $u \in H^2(\Omega)$.}   
\end{eqnarray}
Let  $u_0 \in H^2(\Omega)$ and $\ell(a,b)= \frac{1}{4}\|a\|^2_H-\langle a,u_0 \rangle +\|b-u_0\|_H^2$, then  the following functional 
\begin{equation}
I(u)= \int_0^T \Big [  L  (u(t),\dot { u}(t)+\Lambda  u(t)) +\langle \Lambda  u(t),  u(t) \rangle \Big ] \, dt +\ell (u(0)- u(T), \frac { u(T)+ u(0)}{2})
\end{equation}
attains its minimum at $v \in {\cal X}_{p,q}$ in such a way that $I(v)=\inf_{u\in  {\cal X}_{p,q}}I(u)=0$ and

\begin{eqnarray}\label{principle3}
 \left \{ \begin{array}{lcl}
  \hfill  \dot {v}(t)-i \Delta v(t)+i|v(t)|^{r-1}v(t)& = & -\bar \partial L (v(t)), \\ \label{eqn:zero2}
  \hfill  v(0)&=&u_0.
\end{array}\right.
  \end{eqnarray}
\end{theorem}
\paragraph{Proof} Let $X=H_0^1(\Omega)$ and $H=L^2(\Omega).$ 
Taking into account Theprem \ref{nonlinear}, we just need to verify (\ref{Lambda-10}), (\ref{Lambda-20}) and prove that $(L, \ell)$ is $\Lambda-$coercive on ${\cal X}_{p,q}.$  (\ref{Lambda-20}) follows from  the fact that $\langle \Lambda u, u \rangle =0.$ To prove (\ref{Lambda-10}), note that 
\begin{eqnarray*}
\|\Lambda u\|_{H^{-1}}=\|-\Delta u+|u|^{r-1}u \|_{H^{-1}}\leq \|-\Delta u\|_{H^{-1}}+C\||u|^{r-1}u \|_{L^q(\Omega)}=\|u\|_{H_0^1}+C\|u\|^r_{L^{rq}}.
\end{eqnarray*}
Since $p\geq 2$, we have $qr \leq 2r \leq \frac {2N}{N-2}$.  It follows from Sobolev inequality and the above  that 
\[\|\Lambda u\|_{H^{-1}}\leq \|u\|_{H_0^1} +C \|u\|^r_{H_0^1} \]
from which we obtain
\begin{eqnarray*}
\|\Lambda u\|_{L^q_{H^-1}}  \leq  \|u\|_{L^q_{H_0^1}}+ C\|u\|^r_{L^{rq}_{H_0^1}} 
 \leq  C( \|u\|_{L^p_{H_0^1}}+\|u\|^r_{L^p_{H_0^1}}).
\end{eqnarray*}
To show that $(L, \ell)$ is $\Lambda-$coercive, we assume that $u_n$ is a sequence in ${\cal X}_{p,q}$ such that 
\begin{eqnarray}\label{app-equ}
 \left \{ \begin{array}{lcl}
  \hfill - \dot {u}_n(t)+ i \Delta u_n(t)-i|u_n(t)|^{r-1}u_n(t)& = & -\frac{1}{n}\|u_n\|^{p-2}\Delta u_n +\bar \partial L (u_n(t)), \\ \label{eqn:zero2}
  \hfill  u_n(0)&=&u_0.
\end{array}\right.
  \end{eqnarray}
Since $u_0 \in H^2(\Omega)$, it is standard that at least $u_n \in H^2(\Omega).$ Now if multiply both sides of the above equation by $\Delta u_n(t)-|u_n(t)|^{r-1}u_n(t)$ and taking into account (\ref{pos}) we have 
\begin{eqnarray*}
\langle  \dot {u}_n(t), -\Delta u_n(t)+|u_n(t)|^{r-1}u_n(t)\leq 0
\end{eqnarray*}
from which we obtain 
\begin{eqnarray*}
\frac{1}{2}\|u_n(t)\|^2_{H_0^1}+\frac{1}{r+1}\|u_n(t)\|^{r+1}\leq \frac{1}{2}\|u(0)\|^2_{H_0^1}+\frac{1}{r+1}\|u(0)\|^{r+1}
\end{eqnarray*}
which combined  with (\ref{app-equ}), gives  the boundedness of $u_n$ in ${\cal X}_{p,q}.$  

\paragraph{Example 9} Here are two typical examples for anti-selfdual Lagrangians satisfying the assumptions of the above Theorem

\begin{itemize}
\item $L(u,p)=\phi(u)+\phi^*(-p)$ where $\phi=0$ which leads to a solution of:
\begin{eqnarray*}
 \left \{ \begin{array}{lcl}
  \hfill  i\dot {v}(t)+\Delta v(t)+|v(t)|^{r-1}v(t)& = &0, \\ \label{eqn:zero2}
  \hfill  v(0)&=&u_0
\end{array}\right.
  \end{eqnarray*}
\item $L(u,p)=\phi(u)+\phi^*(a.\nabla u-p)$ where $\phi(u)=\frac{1}{2}\int_{\Omega}|\nabla u|^2 \, dx$ and $a$ is a  vector field on $\Omega$ with compact support. In this case we have a solution for 
\begin{eqnarray*}
 \left \{ \begin{array}{lcl}
  \hfill  i\dot {v}(t)+\Delta v(t)+|v(t)|^{r-1}v(t)& = &-ia.\nabla v+i\Delta v(t), \\ \label{eqn:zero2}
  \hfill  v(0)&=&u_0.
\end{array}\right.
  \end{eqnarray*}
\end{itemize}

\subsection{Noncoercive nonlinear evolutions} 
We shall now assume that there is a symmetric linear duality map $D$ between $X$ and $X^*.$

 \begin{theorem}\label{S-nonlinear}
 Let $(\bar S_t)_{t \in \R}$ be a $C_0-$unitary group of operators associated to a skew-adjoint operator $A$ on the Hilbert space $X^*$, and let $(S_t)_{t \in \R}$ be the corresponding group on $X$. For $p>1$ and $q=\frac{p}{p-1}$, assume that $\Lambda:{\cal X}_{p,q}\to  L^q_{X^{*}}$  is  a  regular map such that for some nondecreasing continuous real function $w$, and $0\leq k<1$, it satisfies
 \begin{equation} \label{Lambda-100}
\hbox{$\|\Lambda S_t x\|_{L^q_{X^*}}\leq k\|\dot x\|_{L^q_{X^*}}+ w(\|x\|_{L^p_{X}})$ for every $x\in {\cal X}_{p,q}$,}
\end{equation}
and
\begin{eqnarray}\label{Lambda-200}
\hbox{$\big |\int_0^T \langle \Lambda  x(t),  x(t)\rangle \, dt  \big | \leq  w(\|x\|_{L^p_X})$ for every $x\in {\cal X}_{p,q}$}.
\end{eqnarray}
 Let $\ell$ be an anti-selfdual Lagrangian on $H \times H$ that is bounded below with  $0 \in {\rm Dom}(\ell)$, and  let $L$ be a time dependent  anti-selfdual Lagrangian on  $[0,T] \times X \times X^*$  such  that for some $C>0$ and $r>1$, we have 
\begin{equation}\label{L.C0}
\hbox{$-C \leq \int_0^TL(t, u(t), 0) dt \leq C(1+\|u\|^r_{L^p_X})$ for every $u\in L^p_X$.}
\end{equation} 
Assume that $(L, \ell)$ is  $\Lambda$-coercive, then the functional 
\begin{equation}
I(u)= \int_0^T \Big [  L  (t,  S_tu(t),\bar S_t\dot { u}(t)+\Lambda  S_tu(t)) +\langle \Lambda S_t u(t),  S_tu(t) \rangle \Big ] \, dt +\ell (u(0)- u(T), \frac { u(T)+ u(0)}{2})
\end{equation}
attains its minimum at $u \in {\cal X}_{p,q}$ in such a way that $I(u)=\inf_{w\in  {\cal X}_{p,q}}I(w)=0$.\\
Moreover if $S_t=\bar S_t$ on $X$, then  $v(t)=S_tu(t)$ is a solution of 
\begin{eqnarray}\label{S-principle3}
 \left \{ \begin{array}{lcl}
  \hfill \Lambda v(t)+A v(t)+ \dot {v}(t)& = &- \bar \partial L (t,v(t)), \\ \label{eqn:zero2}
  \hfill  \frac{v(0)+S_{(-T)}v(T)}{2}
 & \in & -\bar \partial \ell\big(v(0)-S_{(-T)}v(T)).
\end{array}\right.
  \end{eqnarray}
\end{theorem}

\paragraph{Proof:} Define  the nonlinear map $ \Gamma: {\cal X}_{p,q}\to  L^q_{X^{*}}$ by $\Gamma (u)=S_t^* \Lambda S_t (u)$.  This map is also regular in view of the  regularity of  $\Lambda.$  It follows from the previous Lemma that the anti-selfdual Lagrangian $L_{S}$  satisfies (\ref{L.C}). It remains to show that $\Gamma$ satisfies  condition (\ref{Lambda-10}) and (\ref{Lambda-20}). Indeed for $x\in {\cal X}_{p,q}$, we have 
 \begin{equation*} 
\|\Gamma x\|_{L^q_{X^*}}=\|S_t^*\Lambda  S_tx\|_{L^q_{X^*}}=\|\Lambda  S_tx\|_{L^q_{X^*}}\leq k\|\dot x\|_{L^q_{X^*}}+ w(\|x\|_{L^p_{X}})  
\end{equation*}
and
\begin{eqnarray*}
\big |\int_0^T \langle \Gamma  x(t),  x(t)\rangle \, dt  \big |=\big |\int_0^T \langle \Lambda S_t x(t),  S_t x(t)\rangle \, dt  \big | \leq  w(\|S_t x\|_{L^p_X})=w(\|x\|_{L^p_X}). 
\end{eqnarray*}
Also it is easily seen that $L_S$ is $\Gamma-$coercive, which means that 
 all the hypothesis in Theorem \ref{nonlinear} are satisfied. Hence  there exists $u \in {\cal X}_{p,q}$ such that $I(u)=0$ and as in the proof of Theorem \ref{no.diffusion},  $v(t)=S_t u(t)$ is a solution of (\ref{S-principle3}).  

\subsection{Variational resolution for a Fluid driven by $-i\Delta^2$ }

As a consequence of Theorem \ref{nonlinear}, we have provided in \cite{GM3} a variational resolution to evolution equations involving  nonlinear operators such as the Navier-Stokes equation with various boundary conditions. Indeed, by considering  
  \begin{equation}
\label{T-NS.in}
 \left\{ \begin{array}{lcl}
    \hfill
 \frac {\partial u}{\partial t}+(u\cdot \nabla)u +f &=&\nu \Delta u - \nabla  p \quad \hbox{\rm on $ \Omega\subset \R^n$},\\
\hfill {\rm div} \, u&=&0 \quad \hbox{\rm on  $\Omega$},\\
\hfill u&=&0 \quad \hbox{\rm on $\partial \Omega$},\\
\end{array}\right.
\end{equation}
where   $f\in L^{2}_{X^*}([0,T])$,   $X=\{u\in H^{1}_0(\Omega; {\bf R}^{n}); {\rm div} v=0\}$, and $H=L^2(\Omega)$.  Letting 
\begin{equation*}
\label{Phi}
\Phi (u)=\frac{\nu}{2} \int_{\Omega}\Sigma_{j,k=1}^{3}(\frac {\partial u_{j} } {\partial x_{k}})^{2}\, dx+\int_{\Omega}\Sigma_{j=1}^{3}f_{j}u_{j}
\end{equation*}
be the convex continuous function on the space $X=\{u\in H^{1}_0(\Omega; {\bf R}^{n}); {\rm div} v=0\}$, and $\Phi^{*}$ be its Legendre transform  on $X^*$. Equation (\ref{T-NS.in}) can then be reformulated as
 \begin{equation}
\label{T-NS-eq}
 \left\{ \begin{array}{lcl}
  \hfill \frac{\partial u}{\partial t}+ \Lambda u   &\in& -\partial \Phi (t,u) \\
 \hfill  \frac{ u(0)+ u(T)}{2} & \in& -\bar \partial \ell (u(0)-u(T)).
\end{array}\right.
\end{equation}
where  $\Lambda: X \to X^{*}$ is the regular nonlinear operator defined as
\begin{equation}\label{Lambda}
\langle \Lambda u, v\rangle =\int_{\Omega}\Sigma_{j,k=1}^{3}u_{k}\frac {\partial u_{j} } {\partial x_{k}}v_{j}\, dx=\langle (u\cdot \nabla)u,v\rangle.
\end{equation}
and  where $\ell$ is any anti-selfdual Lagrangian on $H \times H$. 
Note that $\Lambda$ maps $X$ into its dual $X^*$ as long as the dimension $N\leq 4$. Moreover, if we lift  $\Lambda$  to path space by defining $(\Lambda u)(t)=\Lambda (u (t))$,  then in dimension $N=2$,  $\Lambda$ is a regular map from ${\cal X}_{2,2}[0,T]$ into  $L^2_{X^*}[0, T]$.
 
 It follows that for  $f$ in $ L^{2}_{X^*}([0,T])$, and if $\psi$ is any convex lower semi-continuous function  on $H$ that is bounded below with $0\in {\rm dom} (\phi)$, then the infimum of the functional
\[
I(u)=\int_0^T \big [\Phi( t, u(t))+\Phi^*(t, -\dot{u}(t)- (u\cdot \nabla)u(t)) \big ] \,dt + \ell (u(0)-u(T), \frac{u(0)+ u(T)}{2})
\]
on ${\cal X}_{2, 2}$ is zero and is attained at a  solution  $u$ of  (\ref{T-NS.in}) that satisfies the following time-boundary condition: 
\begin{equation} \label{br}
 \frac{u(0)+u(T)}{2}
  \in -\bar \partial \ell\big(u(0)-u(T)).
  \end{equation}
   Moreover,  $u$ verifies the following ``energy identity":
\begin{equation}\label{T-NS2-in}
\hbox {$\|u(t)\|_{H}^2+2 \int_0^t \big [\Phi ( t, u(t))+\Phi^*(t, -\dot{u}(t)- (u\cdot \nabla)u(t)) \big ] \,dt = \|u(0)\|_{H}^2$ for every $t \in [0,T].$}
\end{equation}
 Consider now the problem of finding periodic type solutions for the following equation 
 \begin{equation}
\label{T-NS.per}
 \left\{ \begin{array}{lcl}
    \hfill
 \frac {\partial u}{\partial t}+(u\cdot \nabla)u -i\Delta^2 u+f &=&\nu \Delta u - \nabla  p \quad \hbox{\rm on $ \Omega\subset \R^n$},\\
\hfill {\rm div} \, u&=&0 \quad \hbox{\rm on  $\Omega$},\\
\hfill u&=&0 \quad \hbox{\rm on $\partial \Omega$},\\
\end{array}\right.
\end{equation}
where $u=(u_1,u_2)$ and $i\Delta^2 u=(\Delta^2 u_2,-\Delta^2 u_1)$ with $$Dom(i\Delta^2)=\{u \in H_0^1 (\Omega); \Delta u \in H_0^1(\Omega) \text{ and } u=\Delta u=0 \text{ on } \partial \Omega\}.$$
\begin{theorem}\label{T-NS-per} Let $(S_t)_{t \in \R}$ be the  $C_0-$unitary group of operators associtaed to the skew-adjoint operator 
$i\Delta^2.$  Assuming $N=2$, $f$ in $ L^{2}_{X^*}([0,T])$, and $\ell$ to be an anti-selfdual Lagrangian on $H \times H$ that is bounded from below, then the infimum of the functional
\[
I(u)=\int_0^T \big [\Phi( t, S_tu(t))+\Phi^*(t, -S_t\dot{u}(t)- S^*_t\Lambda S_t u(t)) \big ] \,dt + \ell (u(0)-u(T), \frac{u(0)+ u(T)}{2})
\]
on ${\cal X}_{2, 2}$ is zero and is attained at  $u(t)$ in such a way that   $v(t)=S_tu(t)$ is a solution of   (\ref{T-NS.per}) that satisfies the following time-boundary condition: 
\begin{equation} \label{br}
 -\frac{v(0)+S_{(-T)}v(T)}{2}
  \in  \bar \partial \ell\big(v(0)-S_{(-T)}v(T)).
  \end{equation}
   Moreover,  $u$ verifies the following ``energy identity":
\begin{equation}\label{T-NS2-in}
\hbox {$\|u(t)\|_{H}^2+2 \int_0^t \big [\Phi ( t, S_tu(t))+\Phi^*(t, -S_t\dot{u}(t)- S^*_t\Lambda S_t u(t)) \big ] \,dt = \|u(0)\|_{H}^2$ for every $t \in [0,T].$}
\end{equation}
In particular, with appropriate choices for the boundary Lagrangian $\ell$, the solution $v$ can be chosen to  verify either one of the following boundary conditions:
\begin{itemize}
\item an initial value problem:  $v(0) = v_0$ where $v_0$ is a given function in $H$. 
\item a periodic orbit  :  $v(0) = S_{(-T)}v(T)$,
\item an anti-periodic orbit  :  $v(0) =-  S_{(-T)}v(T)$.
 \end{itemize}
\end{theorem}
\paragraph{Proof:} The duality map between $X$ and $X^*$ is $D=- \Delta$ and is therefore linear and symmetric. Also we have $S_t=e^{it\Delta^2}$ and s therefore $S_tD=DS_t.$ Now the result follows from Theorem \ref{S-nonlinear} and the remarks preceeding it.


\begin{thebibliography}{99}
  
  \bibitem{Au1} G. Auchmuty. {\em Saddle points and existence-uniqueness for
evolution equations}, Differential Integral Equations, {\bf 6} (1993),
1161--1171.
\bibitem{Au2} G. Auchmuty. {\em Variational principles for operator equations and initial value problems}, Nonlinear Analysis, Theory, Methods and Applications Vol. {\bf 12}, No.5, pp. 531-564 (1988).

  \bibitem{BE1} H. Brezis, I. Ekeland, {\em Un principe variationnel
  associ\'e \`a certaines equations paraboliques. Le cas independant du
temps}, C.R. Acad. Sci. Paris S\'er. A {\bf 282} (1976), 971--974.

\bibitem{BNS} H. Brezis, L. Nirenberg, G. Stampachia, {\em A remark on Ky Fan's Minimax Principle}, Bollettino U. M. I (1972), 293-300.

\bibitem{Ba} V. Barbu, {\em Abstract periodic Hamiltonian systems},  Adv. Differential Equations  1  (1996),  no. 4, 675--688. 

\bibitem{Br} H. Brezis,  {\em Operateurs maximaux monotones et semi-groupes
de contractions dans les espaces de Hilbert},  North Holland, 
Amsterdam-London, 1973.

\bibitem{Ca} T. Cazenave, {\em Semilinear Schršdinger equations}, Courant Lecture Notes in Mathematics, 10. New York University, Courant Institute of Mathematical Sciences, New York; American Mathematical Society, Providence, RI, 2003. 323 pp.  


  \bibitem{ET} I. Ekeland, R. Temam, {\em Convex Analysis and Variational problems},
Classics in Applied Mathematics, {\bf 28} SIAM (1999 Edition).
  
  \bibitem{Ek.book} I. Ekeland: {\it Convexity Methods in Hamiltonian Mechanics}.
Springer-Verlag, Berlin, Heidelberg, New-York (1990).\par


 \bibitem{G2} N. Ghoussoub,  {\em Anti-selfdual Lagrangians: Variational resolutions of non self-adjoint equations and dissipative evolutions}, AIHP-Analyse non lin\'eaire, 24 (2007) p.171-205.
 
 \bibitem{G3} N. Ghoussoub, {\it Anti-symmetric Hamiltonians: Variational resolution of Navier-Stokes equations and other nonlinear evolutions},  Comm. Pure \& Applied Math., vol. 60, no. 5 (2007) pp. 619-653
 
 \bibitem{G4} N. Ghoussoub, {\it  Selfdual partial differential systems and their variational principles}, Research monograph, In preparation (2006)
 
 \bibitem{G5} N. Ghoussoub, {\it  Maximal monotone operators are selfdual vector fields and vice-versa}, Proc. AMS, in press (2006) 9 pages.

\bibitem{GM1} N. Ghoussoub, A. Moameni, {\it On the existence of Hamiltonian paths connecting Lagrangian submanifolds}, Submitted (2005)

\bibitem{GM2} N. Ghoussoub, A. Moameni, {\it Selfdual variational principles for periodic solutions of Hamiltonian and other dynamical systems}, Comm. in PDE  32, (2007) p. 771-795

\bibitem{GM3} N. Ghoussoub, , A. Moameni, {\it Anti-symmetric Hamiltonians (II): Variational resolution for Navier-Stokes equations and other nonlinear evolutions}, Submitted (2007).

\bibitem{GT1} N. Ghoussoub, L. Tzou. {\em A variational principle for gradient flows}, Math. Annalen, Vol 30, 3 (2004) p. 519-549.

\bibitem{S} R. E. Showalter, { \it Monotone operators in Banach Space and nonlinear partial differential equations}
Math. Surv. Mono. Vol. 49, Am. Math. Soc., Providence,  1997.

\bibitem{St} M. Struwe: {\it  Variational methods and their applications to
non-linear partial differential equations and Hamiltonian systems}.
Springer-Verlag (1990).\par


 \bibitem{Te} R. Temam, {\it Infinite-dimensional dynamical systems in mechanics and physics}, Applied mathematical sciences, 68, Springer-Verlag (1997).

\bibitem{V} I. I. Vrabie, { \it Periodic solutions for nonlinear evolution equations in a Banach space,}
Proc. Amer. Math. Soc. 109 (1990), no. 3, 653-661.

\end{thebibliography}
\end{document}